\documentclass[a4paper,oneside,12pt]{article}

\usepackage{latexsym}
\usepackage[latin1]{inputenc}
\usepackage[sumlimits]{amsmath}
\usepackage{amsthm}
\usepackage{amstext}
\usepackage[dvips]{graphicx}
\usepackage{amssymb}
\usepackage{fancyhdr}

\newtheorem{defn}{Definition}[section]
\newtheorem{theorem}[defn]{Theorem}
\newtheorem{corollario}[defn]{Corollary}
\newtheorem{lemma}[defn]{Lemma}
\newtheorem{proposition}[defn]{Proposition}

\newtheorem{oss}[defn]{Remark}

\def\Rset{\mathbb{R}}
\def\Cset{\mathbb{C}}
\def\Nset{\mathbb{N}}
\def\Hset{\mathbb{H}}
\def\Bset{\mathbb{B}}

\newcommand{\abs}[1]{\left|#1\right|}	
\newcommand{\decl}{:=}			

\newcommand{\be}{\begin{equation}}	
\newcommand{\ee}{\end{equation}}	
\newcommand{\bt}{\begin{theorem}}	
\newcommand{\et}{\end{theorem}} 	
\newcommand{\bd}{\begin{defn}}	
\newcommand{\ed}{\end{defn}}	
\newcommand{\bc}{\begin{corollario}}	
\newcommand{\ec}{\end{corollario}}	
\newcommand{\bp}{\begin{dimostrazione}}	
\newcommand{\ep}{\end{dimostrazione}}	
\newcommand{\bl}{\begin{lemma}}		
\newcommand{\el}{\end{lemma}}		
\newcommand{\bpr}{\begin{proposition}}	
\newcommand{\epr}{\end{proposition}}	
\newcommand{\bern}{\begin{eqnarray*}}	
\newcommand{\eern}{\end{eqnarray*}}	
\newcommand{\boss}{\begin{oss}\rm}	
\newcommand{\eoss}{\end{oss}}	

\newcommand{\besi}[1][{}]
	{\begin{esempi}\rm\textbf{#1}\begin{enumerate}}	
\newcommand{\eesi}{\end{enumerate}\end{esempi}}	
\newcommand{\beso[1]}{\begin{esempio}[#1]\rm}	
\newcommand{\eeso}{\end{esempio}}

\newcommand{\Cuno}{\mbox{\cal C}^{\,1}}
\newcommand{\C}{\mbox{\cal C}}
\newcommand{\Cinfinito}{\mbox{\cal C}^{\,\infty}}

\renewcommand{\theequation}{\thesection.\arabic{equation}}

\newcounter{proofeqno}		
\newenvironment{dimostrazione}	
{\noindent \textbf{Proof.\ }\setcounter{proofeqno}{0}}
{\hfill $\Box$\par\medskip}

\def\acknowledgment{\goodbreak\subsection*{Acknowledgment}
\bgroup \footnotesize}

\def\acknowledgments{\goodbreak\subsection*{Acknowledgments}
\bgroup \footnotesize}

\def\endacknowledgment{\vskip1sp\egroup}
\def\endacknowledgments{\vskip1sp\egroup}

\newcommand {\Div} {\mbox{\rm div\,}}

\newcommand{\appendice}[2]{
\renewcommand{\theequation}{#2-\arabic{equation}}
\setcounter{equation}{0}
\setcounter{section}{0}
\setcounter{defn}{0}
\renewcommand{\thesection}{#2-\arabic{section}}
\section*{#2  #1}
\renewcommand{\thedefn}{#2-\arabic{defn}}
\addcontentsline{toc}{section}{#2 #1} {}
\markboth{#1}{}
}

\title{Hardy inequalities on Riemannian manifolds \\ and applications }
\author{Lorenzo D'Ambrosio\\{\small Dipartimento di Matematica, via E. Orabona, 4 -- I-70125, Bari, Italy}\thanks{Address for correspondence e-mail: dambros@dm.uniba.it}\\ \\
Serena Dipierro\\{\small SISSA, Sector of Mathematical Analysis, via Bonomea 265 -- I-34136, Trieste, Italy}\thanks{Address for correspondence e-mail: serydipierro@yahoo.it}}

\date{April 6, 2012}

\begin{document}
\maketitle
{\abstract{We prove a simple sufficient criteria to obtain some Hardy inequalities on 
Riemannian manifolds related to quasilinear 
second-order differential operator $\Delta_{p}u := \Div\left(\abs{\nabla u}^{p-2}\nabla u\right)$. 
Namely, if $\rho$ is a nonnegative weight such that $-\Delta_{p}\rho\geq0$, 
then the Hardy inequality
\be c\int_{M}\frac{\abs{u}^{p}}{\rho^{p}}\abs{\nabla \rho}^{p} dv_{g} \leq
\int_{M}\abs{\nabla u}^{p} dv_{g}, \quad u\in\Cinfinito_{0}(M).
\nonumber\ee 
holds. We show concrete examples specializing the function $\rho$.

Our approach allows to obtain a characterization of $p$-hyperbolic manifolds
as well as other inequalities related to  Caccioppoli inequalities, weighted Gagliardo-Nirenberg inequalities, uncertain principle and first order Caffarelli-Kohn-Nirenberg interpolation inequality.

}

}

\bigskip



\medskip
\tableofcontents

\section{Introduction}

An $N$-dimensional generalization of the classical 
Hardy inequality asserts that for every $p>1$
\be 
c \int_{\Omega}\frac{\abs{u}^{p}}{w^{p}} \, dx \leq\int_{\Omega}\abs{\nabla u}^{p} dx, 
\quad   u\in\Cinfinito_{0}\left(\Omega\right), \nonumber\ee
where $\Omega\subset\Rset^{N}$ is an open set, 
and the weight $w$ is, for instance, 
$w\left(x\right) \decl \abs{x}$ and $p<N$, 
or $w\left(x\right) \decl dist\left(x,\partial\Omega\right)$ and $\Omega$ is convex 
(see for example \cite{BFT, BM, DH, GGM, MMP, MS1, MS2, Mi} and references therein). 

The preeminent role of Hardy inequalities and the knowledge of the best constants 
involved is a well known fact, as the reader can recognize from the wide literature 
that uses such a tool in Euclidean or in subelliptic setting as well as on manifolds 
(\cite{BG, BGGK, BC, BM, BD, DL, Ja, MP} just to cite a few).

On the other hand, 
the knowledge of the validity of a Hardy or Gagliardo-Nirenberg or Sobolev or Caffarelli-Kohn-Nirenberg inequality on a manifold $M$
and their best constants allows to obtain qualitative properties on the manifold $M$. 
For instance in \cite{AX,CX,Xi} it was shown that if $M$ is a complete open Riemannian manifold 
with nonnegative Ricci curvature in which a Hardy or Gagliardo-Nirenberg 
or Caffarelli-Kohn-Nirenberg type inequality holds, 
then $M$ is in some suitable sense \emph{close} to the Euclidean space. 

\bigskip

One of our aims is to prove some Hardy inequalities on Riemannian manifolds. 
In 1997, Carron in \cite{Ca} studies weighted $L^{2}$-Hardy inequalities 
on a Riemannian manifold $M$
under some geometric assumptions on the weight function $\rho$, 
obtaining, among other results, the following inequality 
\be  
c \int_{M} \frac{u^{2}}{\rho^{2}}\, dv_{g}\leq\int_{M}\abs{\nabla u}^{2} dv_{g}, 
\quad u\in\Cinfinito_{0}\left(M\right), \label{carron}\ee
where $\rho$ is a nonnegative function such that $\abs{\nabla\rho}=1$, 
$\Delta\rho\geq \frac{\gamma}{\rho}$, $\rho^{-1}\left\{0\right\}$ 
is a compact set of zero capacity and 
$c=\left(\frac{\gamma -1}{2}\right)^{2}$. 
In \cite{Ca} the author applies this result to several explicit examples 
of Riemannian manifolds. 
Under the same hypotheses on the function~$\rho$, Kombe and \"Ozaydin
in \cite{KO} extend 
Carron's result to the case~$p\neq2$ for functions in~$\Cinfinito_{0}\left(M\setminus\rho^{-1}\left\{0\right\}\right)$, 
and the authors present an application to the punctured 
manifold $\Bset^n\setminus\left\lbrace0\right\rbrace$ 
with $\Bset^n$ the Poincar\'e ball model of the hyperbolic space
and $\rho$ the distance from the point $0$ and $p=2$.

Li and Wang in \cite{LW} prove that if $M$ is a hyperbolic manifold 
(i.e. there exists a symmetric positive Green function $G_{x}\left(\cdot\right)$ 
for the Laplacian with pole at $x$), then
\be \frac{1}{4} \int_{M} \frac{\abs{\nabla_{y} G_{x}\left(y\right)}^{2}}{G_{x}^{2}\left(y\right)} \, u^{2}(y) \, dv_{g} \leq\int_{M}\abs{\nabla u}^{2}(y) \, dv_{g}, \quad u\in\Cinfinito_{0}\left(M\setminus\left\{x\right\}\right). \nonumber \ee 

We also mention Miklyukov and Vuorinen, which in \cite{MV} prove that 
the inequality
\be \left(\int_{M} \abs{\alpha\left(\epsilon\left(x\right)\right)u\left(x\right)}^{q} dv_{g}\right)^{1/q} \leq \lambda \left(\int_{M} \left(\beta\left(\epsilon\left(x\right)\right)\abs{\nabla u\left(x\right)}\right)^{p} dv_{g}\right)^{1/p}, \quad u\in W^{1,p}_{0}\left(M\right),\nonumber\ee
holds for $q\geq p$ provided  some conditions related to the isoperimetric profile of $M$ are satisfied. 

In \cite{AS}, Adimurthi and Sekar use the fundamental solution of a general 
second-order elliptic operator to derive Hardy-type inequalities and then they 
extend their arguments to Riemannian manifolds using the fundamental 
solution of $p$-Laplacian. 

Bozhkov and Mitidieri in \cite{BMi1} prove the validity of $(\ref{carron})$ also for $p\neq 2$ 
($1<p<N$), provided there exists on $M$ a $\Cuno$ conformal Killing vector field 
$K$ such that $\Div K=\mu$ with $\mu$ a positive constant and $\rho=\abs{K}$. 

Let $p>1$ and let $\rho$ be a nonnegative function. 
Our principal result is a simple criterion to establish if there holds 
a Hardy inequality involving the weight $\rho$.
Namely, if $\rho$ is $p$-superharmonic in $\Omega$,  
that is $-\Delta_{p}\rho\geq0$, then
\be c \int_{M}\frac{\abs{u}^{p}}{\rho^{p}}\abs{\nabla \rho}^{p} dv_{g}\leq \int_{M}\abs{\nabla u}^{p} dv_{g}, \quad u\in\Cinfinito_{0}\left(M\right), \label{d5}\ee
holds (see Theorem $\ref{teo:g}$). 
Such a kind of criteria is already established in \cite{Da4} for 
a quite general class of second order operators containing, among other examples, 
the subelliptic operators on Carnot groups. 
For this goal we shall mainly use a technique introduced by Mitidieri in \cite{Mi} 
and developed in \cite{Da2, Da3, Da4} and in \cite{BMi1, BMi2}. 
The proof is based on the divergence theorem and on the careful choice of a vector field. 

\medskip

Let us point out some interesting outcomes of our approach. 
A first issue is that, since it is quite general, 
our approach includes Hardy inequalities already 
studied in \cite{AS, BMi1, Ca, KO, LW} in the case $p=2$ 
as well as their generalization for $p>1$. Indeed, in all these cited papers,
the authors assume extra conditions on the function $\rho$ or on the manifold. 
Furthermore, in concrete cases, 
our result yields an explicit value of the constant $c$. 
Moreover, in several cases, this value is also the best constant (see \cite{Da4}). 
To this regard, we discuss if the best constant is achieved or not and, in the latter case, we study the possibility to add a remainder term. 

Another aspect of our technique is that it allows to characterize the 
$p$-hyperbolic manifolds. We remind that
a manifold $M$ is called \emph{$p$-hyperbolic} if there exists a symmetric 
positive Green function $G_{x}\left(\cdot\right)$ for 
the $p$-Laplacian with pole at $x$. 
We prove that $M$ is $p$-hyperbolic if and only if there exists a nonnegative
non trivial function $f\in L^1_{loc}\left(M\right)$ such that
\[\int_{M} f\, {\abs{u}^{p}} dv_{g}\leq \int_{M}\abs{\nabla u}^{p} dv_{g}, \quad u\in\Cinfinito_{0}\left(M\right). 
\]
Notice that one of the implications of this characterization for $p=2$ is 
the result proved in \cite{LW}. 
During the review process of this work, 
we have received the paper of Devyver, Fraas and Pinchover  \cite{DFP}.
In \cite{DFP} a general linear second order differential 
operator $P$ in the Euclidean framework is studied. The authors
find a profound relation between the existence of positive
supersolutions of $Pu=0$,  Hardy type inequalities involving $P$ and 
a weight $W$ and the characterization of the spectrum of the weighted operator.
We refer the interested reader to \cite{DFP, DFP2}.

We also obtain a generalization of (\ref{d5}). Namely, for a nonnegative function $\rho$, the inequality
\be  \left( \frac{\abs{p-1-\alpha}}{p}\right)^{p} \int_{M} \frac{\abs{u}^{p}}{\rho^{p}}\rho^\alpha \abs{\nabla \rho}^{p} dv_{g} \leq \int_{M}\abs{\nabla u}^{p}\rho^\alpha \, dv_{g}, \quad u\in\Cinfinito_{0}\left(M\right),
\label{gener} \ee
holds, provided $-(p-1-\alpha)\Delta_p \rho \ge 0$ (see Theorem \ref{hardygen}). 
The above inequality contains, as special case,
the Caccioppoli inequality. Indeed, if $\rho$ is a $p$-subharmonic function,
that is $\Delta_p \rho\ge 0$, then $(\ref{gener})$ holds for $\alpha>p-1$ and, 
in particular, for $\alpha=p$ we have 
\[ \frac{1}{p^{p}}  \int_{M} {\abs{u}^{p}} \abs{\nabla \rho}^{p} dv_{g}\leq \int_{M}\abs{\nabla u}^{p}\rho^p \, dv_{g}, \quad u\in\Cinfinito_{0}\left(M\right). 
\]
This is the so called Caccioppoli inequality (see for instance \cite{PRS} and 
the references therein for the version $p=2$ on manifolds). 

Another advantage of our approach is that it allows to obtain also 
other new and known results, like 
wighted Gagliardo-Nirenberg inequalities and the uncertain principle. 

Finally we show that if $(\ref{gener})$ and 
a Sobolev type inequality (that is $c\abs u_{L^{p^*}}\le \abs{\nabla u}_{L^p}$) hold on $M$, 
then we obtain an interpolation inequality involving, as weights,
$\rho$ and its gradient. 
As particular case, our results contain inequalities on manifolds 
related to the celebrated Caffarelli-Kohn-Nirenberg inequality. 

\bigskip
 
The paper is organized as follows. 
We present the proof of $(\ref{d5})$ in Section $\ref{sec:hi}$, where important 
consequences and observations are derived. 
In Section $\ref{sec:gen}$ we show natural extensions of $(\ref{d5})$, 
obtaining also Hardy inequality with weights, 
Caccioppoli-type inequalities, 
weighted Gagliardo-Nirenberg inequalities and the uncertain principle. 
Some remarks on the best constant and if it is attained are discussed in 
Section $\ref{sec:bc}$. 
In Section $\ref{sec:ckn}$ we present a first order interpolation inequality. 
Finally Section $\ref{sec:ex}$ is devoted to present some concrete examples of
Hardy-type  inequalities on manifolds. 

\subsection*{Notation} 
In what follows $\left(M,g\right)$ is 
a complete Riemannian $N$-dimensional manifold, 
$\Omega\subset M$ is an open set, 
$dv_{g}$ is the volume form associated to the metric $g$, 
$\nabla u$ and $\Div h$ stand respectively for the gradient of a function $u$ 
and the divergence of a vector field $h$ 
with respect to the metric $g$ 
(see \cite{Au} for further details). 
Throughout this paper $p>1$. 

\section{Hardy  inequalities}\label{sec:hi}

In order to state Hardy inequalities involving a weight $\rho$, 
the basic assumption we made on $\rho$ is that 
$\rho$ is $p$-superharmonic in weak sense. 
Namely, we assume that $\rho\in L^{1}_{loc}\left(\Omega\right)$, 
$\abs{\nabla\rho}\in L^{p-1}_{loc}\left(\Omega\right)$, 
and $-\Delta_{p}\rho\geq 0$ on $\Omega$ in weak sense, 
that is for every nonnegative $\varphi\in\Cuno_{0}(\Omega)$, we have 
\be\int_{\Omega}\abs{\nabla\rho}^{p-2} \left( \nabla\rho\cdot\nabla\varphi \right)  dv_{g} \geq 0. \label{subarmonica10}\ee

The main result on Hardy inequalities is the following:
\bt \label{teo:g} 
Let $\rho\in W^{1, p}_{loc}\left(\Omega\right)$ be a nonnegative function on $\Omega$ such that 
\be -\Delta_{p}\rho\geq 0 \quad \mathrm{on\ }  \Omega \mathrm{\ in\ weak\ sense} .  \nonumber\ee
Then $\frac{\abs{\nabla\rho}^{p}}{\rho^{p}}\in L^{1}_{loc}\left(\Omega\right)$, 
and the following inequality holds:
\be \left(\frac{p-1}{p}\right)^{p}\int_{\Omega}\frac{\abs{u}^{p}}{\rho^{p}}\abs{\nabla \rho}^{p} dv_{g}\leq \int_{\Omega}\abs{\nabla u}^{p} dv_{g}, 
\quad u\in\Cinfinito_{0}\left(\Omega\right). 
\label{dishardy10}\ee
\et

Before proving Theorem $\ref{teo:g}$, we shall present some immediate consequences 
and extensions of the main result. 

\bd\label{def:a}
Let $\Omega\subset M$ be an open set. 
We denote by $D^{1,p}\left(\Omega\right)$ the completion of 
$\Cinfinito_{0}\left(\Omega\right)$ with respect to the norm
\be \abs{u}_{D^{1,p}} =\left(\int_{\Omega}\abs{\nabla u}^{p}dv_{g}\right)^{1/p}. \nonumber\ee
\ed

It is possible to extend the validity of $(\ref{dishardy10})$ to function $u\in\Cinfinito_{0}\left(M\right)$. 
This extension is based on the inclusion 
\be D^{1,p}\left(M\right) \subset D^{1,p}\left(\Omega\right).   \label{inclus1}\ee
The above inclusion is satisfied, for instance, when $M\setminus\Omega$ is a compact set 
of zero $p$-capacity (see Appendix A). 

\bc\label{coro:q1} 
Let $\rho\in L^{1}_{loc}\left(M\right)$ be a function satisfying the assumptions of Theorem $\ref{teo:g}$. 
If $(\ref{inclus1})$ holds, then $\frac{\abs{\nabla\rho}^{p}}{\rho^{p}}\in L^{1}_{loc}\left(M\right)$, 
and the following inequality holds:
\be\left(\frac{p-1}{p}\right)^{p} \int_{M}\frac{\abs{u}^{p}}{\rho^{p}}\abs{\nabla \rho}^{p} dv_{g} 
\leq \int_{M}\abs{\nabla u}^{p} dv_{g}, \quad u\in\Cinfinito_{0}\left(M\right). \label{dishardy11}\ee
\ec

\noindent{\bf Proof.} 
The inequality $(\ref{dishardy10})$ holds for every $u\in\Cinfinito_{0}\left(\Omega\right)$, 
then it holds for every $u\in D^{1,p}\left(\Omega\right)$. 
Since $\Cinfinito_{0}\left(M\right) \subset  D^{1,p}\left(M\right)$,  
by using $(\ref{inclus1})$ we conclude the proof.
\hfill$\Box$

\bigskip

In order to illustrate further consequences of Theorem $\ref{teo:g}$ we give the following:
\bd\label{def:p-par} 
A manifold $M$ is said \emph{$p$-hyperbolic}\footnote{Many authors call these manifolds non p-parabolic.} if there exists a symmetric positive Green function $G_{x}\left(\cdot\right)$ for the $p$-Laplacian with pole at $x$\footnote{That is $-\Delta_p G_x=\delta_x$ where $\delta_x$ is the Dirac measure concentrated at point $x$. }, if it is not the case we call it \emph{$p$-parabolic}. 
\ed 

Several equivalent definitions of $p$-parabolic manifolds can be given. 
For instance in \cite{Tr1} there is the following (see also the literature therein and \cite{Ho1})
\bpr\label{pr:par} Let $p>1$. The following statements are equivalent
\begin{itemize}
\item[a)] $M$ is $p$-parabolic;
\item[b)] there exists a compact set $K\subset M$ 
with non empty interior such that $cap_p (K, M) = 0$;
\item[c)] there is no non constant positive $p$-superharmonic function on $M$;
\item[d)] there exists a sequence of functions $u_j\in\C^\infty_0(M)$ 
such that $0 \le u_j \le 1$, $u_j \to 1$ uniformly on every compact subset of $M$ and
$\int_M\abs{\nabla u_j}^p dv_g\to 0$.
\end{itemize}
\epr

Other characterizations of $p$-parabolic manifolds are based on several properties, for instance on the volume growth, on the isoperimetric profile of the manifold, on some properties of some cohomology, on the recurrence of the Brownian motion. 
See \cite{Gr1,Gr2,Gr3, Gr4,LT, Tr1} and the references therein.

\medskip

From Theorem \ref{teo:g} we deduce the following characterization of $p$-hyperbolicity.
\bt\label{teo:parabolic}
A manifold $M$ is  $p$-hyperbolic if and only if there exists a nonnegative 
non trivial function $f\in L^1_{loc}\left(M\right)$ such that
\be \int_{M} f\, {\abs{u}^{p}} \, dv_{g}\leq \int_{M}\abs{\nabla u}^{p} dv_{g}, \quad u\in\Cinfinito_{0}\left(M\right). 
\label{par}\ee
\et

\noindent{\bf Proof.} 
If $M$ is $p$-hyperbolic then inequality $(\ref{par})$ holds 
with $f = \left(\frac{p-1}{p}\right)^{p}\frac{\abs{\nabla G_{x}}^{p}}{G_{x}^{p}}$. 
Indeed $G_{x}$ is nonnegative and satisfies the hypotheses of Theorem $\ref{teo:g}$ 
(see Theorem $\ref{teo-par}$ for further details). 

Conversely, assume that $M$ is $p$-parabolic and that $(\ref{par})$ is valid for a function $f\geq0$. 
Then from $d)$ of Proposition \ref{pr:par} 
there exists a sequence of functions 
$u_{j}\in\Cuno_{0}\left(M\right)$ such that $0\leq u_{j}\leq1$, $u_{j}\rightarrow1$ 
uniformly on every compact subset of $M$ and
\be \int_{M}\abs{\nabla u_{j}}^{p} dv_{g} \longrightarrow0, \qquad \mathrm{(as\ } j\rightarrow+\infty).\nonumber\ee 
It implies that $\int_{D}f dv_{g} = 0$ for every compact subset $D$ of $M$ and then 
$f\equiv0$. This concludes the proof.
\hfill$\Box$

\boss
Since the $p$-hyperbolicity of $M$ is equivalent to the existence of a non constant positive 
$p$-superharmonic function $\rho$ on $M$, then by Theorem $\ref{teo:g}$ we obtain that inequality $(\ref{par})$ holds 
with $f = \left(\frac{p-1}{p}\right)^{p}\frac{\abs{\nabla\rho}^{p}}{\rho^{p}}$. 
\eoss

\boss 
Our Theorem \ref{teo:parabolic} implies that if the manifold $M$ admits a $\Cuno$ conformal Killing vector field $K$ (see i.e. \cite{BMi1} for the definition) 
such that $\Div K=\mu\neq 0$ with $\mu$ constant and $\abs K^{-p}\in L^1_{loc}(M)$, 
then $M$ is $p$-hyperbolic. 
This follows combining Theorem $\ref{teo:parabolic}$ and Theorem $4$ of \cite{BMi1} (see also Remark \ref{remspecial} ii) below).
\eoss

\bigskip

In order to prove Theorem $\ref{teo:g}$ we fix some notation. 
Let $h\in L^{1}_{loc}\left(\Omega\right)$ be a vector field. 
We remind that the distribution $\Div h$ is defined as  
\be \int_{\Omega} \varphi \ \Div h \ dv_{g} = -\int_{\Omega}  \left( \nabla{\varphi}\cdot h\right)  dv_{g}, \label{div10}\ee
for every $\varphi\in\Cuno_{0}\left(\Omega\right)$.

Let $h\in L^{1}_{loc}\left(\Omega\right)$ be a vector field and 
let $A\in L^{1}_{loc}\left(\Omega\right)$ be a function.
In what follows we write $A \leq\Div h$ meaning that the inequality holds
in distributional sense, that is for every $\varphi\in\Cuno_{0}\left(\Omega\right)$ such that $\varphi\geq0$, 
we have
\be \int_{\Omega} \varphi \ A \ dv_{g} \leq \int_{\Omega} \varphi \ \Div h \ dv_{g} = -\int_{\Omega} \left( \nabla{\varphi}\cdot h\right) dv_{g}. \label{distr10}\ee

\boss \label{oss:s} 
Let $f\in\Cuno\left(\Rset\right)$ be a real function such that $f\left(0\right) = 0$. 
Taking $\varphi = f\left(u\right)$ with $u\in\Cuno_{0}\left(\Omega\right)$ in (\ref{div10}), 
we have 
\be 
\int_{\Omega} f\left(u\right) \Div h \ dv_{g} = - \int_{\Omega} f'\left(u\right) \left( \nabla{u}\cdot h\right)  dv_{g}.  \nonumber\ee
In particular, choosing $f\left(u\right) = \abs{u}^{p}$, with $p>1$, we get
\be \int_{\Omega} \abs{u}^{p} \Div h \ dv_{g} = -p \int_{\Omega}\abs{u}^{p-2}u \left( \nabla{u}\cdot h\right)  dv_{g}, \quad u\in\Cuno_{0}\left(\Omega\right).  \label{div11}\ee
\eoss

\bl\label{lemma:k1} 
Let $h\in L^{1}_{loc}\left(\Omega\right)$ be a vector field and let 
$A_{h}\in L^{1}_{loc}\left(\Omega\right)$ be a nonnegative function such that
\begin{enumerate} 
\item[i)] $A_{h} \leq\Div h,$
\item[ii)] $\frac{\abs{h}^{p}}{A_{h}^{p-1}}\in L^{1}_{loc}\left(\Omega\right)$. 
\end{enumerate}
Then for every $u\in\Cuno_{0}\left(\Omega\right)$ we have
\be\int_{\Omega} \abs{u}^{p} A_{h} \ dv_{g}\leq p^{p} \int_{\Omega} \frac{\abs{h}^{p}}{A_{h}^{p-1}} \abs{\nabla{u}}^{p} dv_{g}. \label{lem10}\ee
\el

\noindent{\bf Proof.}
We note that the right hand side of $(\ref{lem10})$ is finite 
since $u\in\Cuno_{0}\left(\Omega\right)$. 
Using the identity (\ref{div11}) and H\"older inequality we obtain
\begin{eqnarray} 
\int_{\Omega}\abs{u}^{p} A_{h} \ dv_{g} &\leq& 
\int_{\Omega} \abs{u}^{p} \Div h \ dv_{g} \nonumber\\
&\leq& p \int_{\Omega} \abs{u}^{p-1} \abs{h} \abs{\nabla{u}} dv_{g}  \nonumber\\
&=& p \int_{\Omega} \abs{u}^{p-1} A^{\left(p-1\right)/p}_{h} \frac{\abs{h}}{A^{\left(p-1\right)/p}_{h}} \abs{\nabla{u}} dv_{g} \nonumber\\
&\leq & p \left(\int_{\Omega} \abs{u}^{p} A_{h} \ dv_{g}\right)^{\left(p-1\right)/p}\left(\int_{\Omega} \frac{\abs{h}^{p}}{A_{h}^{p-1}}        \abs{\nabla{u}}^{p} dv_{g}\right)^{1/p}. \nonumber 
\end{eqnarray}
This completes the proof.
\hfill$\Box$

\bigskip

Specializing the vector field $h$ and the function $A_{h}$, 
we shall deduce from $(\ref{lem10})$ Hardy-type inequalities 
on Riemannian manifolds.

\boss\label{remspecial}
Letting us to point out a strategy to get Hardy inequalities 
at least in some special cases. 
Under the hypotheses of Lemma $\ref{lemma:k1}$, if $A_{h} = \Div h$, 
then $(\ref{lem10})$ reads as
\be\int_{\Omega} \abs{u}^{p} \Div h\  dv_{g}\leq p^{p} \int_{\Omega} \frac{\abs{h}^{p}}{\abs{\Div h}^{p-1}} \abs{\nabla{u}}^{p} dv_{g} . \label{lem11}\ee
\begin{itemize}
\item[i)] Let $V$ be a function in $L^{1}_{loc}\left(\Omega\right)$ such that 
its weak partial derivatives of order up to two are in $L^{1}_{loc}\left(\Omega\right)$. 
If $\Delta V\geq 0$, choosing $h = \nabla V$, 
we obtain $A_{h} = \Div h = \Div\left(\nabla V\right) = \Delta V\geq 0$. 
Then from $(\ref{lem11})$
\be \int_{\Omega} \abs{u}^{p} \abs{\Delta V} dv_{g} \leq p^{p} \int_{\Omega} \frac{\abs{\nabla V}^{p}}{\abs{\Delta V}^{p-1}} \abs{\nabla{u}}^{p} dv_{g}. \label{lem12}
\ee
This kind of inequalities for the Euclidean setting $\Omega = \Rset^{N}$ 
are already found by Davies e Hinz in \cite{DH}. 

\noindent At this point, in order to deduce from $(\ref{lem12})$ an inequality like
\be c \int_{\Omega}\frac{\abs{u}^{p}}{\rho^{p}} \, dv_{g} \leq\int_{\Omega}\abs{\nabla u}^{p} dv_{g}, \label{hoss}\ee
we have to choose a suitable function $V$. 
Let us to consider the case when $\rho$ is the 
distance from a point $o\in M$ and $\abs{\nabla\rho}=1$. 
A suitable choice for $V$ is $V = \rho^{2-p}$ if $1<p<2$, 
$V = \ln\rho$ if $p=2$ and $V = -\rho^{2-p}$ if $2<p<N$. 
Following \cite{Ca}, if we require that $\Delta\rho\geq\frac{\gamma}{\rho}$ 
the above choices yield the inequality (\ref{hoss}) with $c = \left(\frac{\gamma-p+1}{p}\right)^{p}$. 
The success of this strategy is deeply linked to the hypothesis $\abs{\nabla\rho}=1$. 
Indeed, it seems that such a strategy does not work even in the subelliptic setting, 
where the analogous of the hypothesis $\abs{\nabla\rho}=1$ does not hold. \\
Furthermore, the fact that the hypothesis $\abs{\nabla\rho}=1$ 
is sometimes restrictive even in the Euclidean case, can be seen in the following example.
In the Euclidean unit ball $B_1\subset\Rset^N$ the inequality
\be c \int_{B_1}\frac{\abs{u}^{p}}{\abs{\abs x \ln \abs x}^p} \, dx \leq\int_{B_1}\abs{\nabla u}^{p} dx \label{hosslog}\ee
holds for $1<p\le N$ (see Section \ref{sec:Cartan}, Section \ref{sec:euclidean} and \cite{Da4}). If we wish to deduce (\ref{hosslog}) from (\ref{hoss})
we are forced to choose $\rho=-\abs x \ln \abs x$. However $\abs{\nabla\rho}\not =1$.
\item[ii)] Let $p<N$. Assume that there exists a $\Cuno$ conformal Killing vector field $K$ (see i.e. \cite{BMi1} for the definition) such that $\Div K=\frac{N}{2}\mu >0$. 
Choosing $h=\frac{K}{\abs{K}^{p}}$, we have $A_{h}:= \Div h = \frac{N-p}{2}\frac{\mu}{\abs{K}^{p}}$ (see Lemma $3$ in \cite{BMi1}) and the inequality $(\ref{lem11})$ reads as 
\be\left(\frac{N-p}{Np}\right)^{p}\int_{\Omega} \frac{\Div K}{\abs{K}^{p}}\abs{u}^{p}  dv_{g} \leq \int_{\Omega} \left(\Div K\right)^{1-p} \abs{\nabla{u}}^{p} dv_{g}. \label{bozmit}\ee 
Therefore, by Lemma $\ref{lemma:k1}$,  $(\ref{bozmit})$ holds for every $u\in\Cuno_{0}\left(\Omega\right)$ provided $\abs{K}^{-p}\in L^{1}_{loc}\left(\Omega\right)$. This last fact was obtained in \cite[Theorem $4$]{BMi1}.
\end{itemize}
\eoss

\noindent{\bf Proof of Theorem \ref{teo:g}.} 
Let $0<\delta <1$, and $\rho_{\delta}:= \rho + \delta$. 
In order to apply Lemma $\ref{lemma:k1}$, we define $h$ and $A_{h}$ as 
\be  h:=  -\frac{\abs{\nabla\rho_\delta}^{p-2}\nabla\rho_\delta}{\rho_{\delta}^{p-1}} 
\quad \mathrm{and\ } A_{h} := \left(p-1\right)\frac{\abs{\nabla\rho_\delta}^{p}}{\rho_{\delta}^{p}}. \label{hAh}\ee
Since $\frac{1}{\rho_\delta}\leq\frac{1}{\delta}$, 
the fact that $\rho\in W^{1, p}_{loc}\left(\Omega\right)$ implies that 
$h\in L^{1}_{loc}\left(\Omega\right)$ 
and $A_{h}\in L^{1}_{loc}\left(\Omega\right)$.
Moreover, by computation we have
\be \frac{\abs{h}^{p}}{A_{h}^{p-1}} = \frac{\abs{\nabla\rho_\delta}^{p\left(p-1\right)}}{{\rho_\delta}^{p\left(p-1\right)}}\frac{{\rho_\delta}^{p\left(p-1\right)}}{\abs{\nabla\rho_\delta}^{p\left(p-1\right)}}\frac{1}{\left(p-1\right)^{\left(p-1\right)}} = \frac{1}{\left(p-1\right)^{p-1}} \in L^{1}_{loc}\left(\Omega\right), \nonumber\ee
that is $ii)$ of Lemma $\ref{lemma:k1}$ is fulfilled. 

The hypothesis $i)$ of Lemma $\ref{lemma:k1}$ is satisfied provided   
\be \left(p-1\right) \int_{\Omega}\frac{\abs{\nabla\rho_\delta}^{p}}{\rho_{\delta}^{p}}\, \varphi \ dv_{g}\leq\int_{\Omega}  \left( \frac{\abs{\nabla\rho_\delta}^{p-2}\nabla\rho_\delta}{\rho_{\delta}^{p-1}}\cdot\nabla\varphi \right) dv_{g}
\label{subarmonica11}\ee
holds for every nonnegative function $\varphi\in\Cuno_{0}\left(\Omega\right)$. 
Then, for a fixed $\varphi\in\Cuno_{0}\left(\Omega\right)$ nonnegative, 
we have to prove $(\ref{subarmonica11})$. 
Let $K = supp\varphi\subset\Omega$ and let $U$ be a neighborhood of $K$ with compact closure in $\Omega$.
We note that both integrals in $(\ref{subarmonica11})$ are finite 
since $\frac{1}{\rho_\delta}\leq\frac{1}{\delta}$ and $\rho\in W^{1, p}_{loc}\left(\Omega\right)$. Since
\be \vert\nabla\ln\rho_{\delta}\vert = \frac{\vert\nabla\rho_{\delta}\vert}{\rho_{\delta}} \leq\frac{\vert\nabla\rho\vert}{\delta}  \in L^{p}_{loc}\left(\Omega\right),
\ee
and $\ln\rho_{\delta}\in L^{p}_{loc}\left(\Omega\right)$, 
we have that $\ln\rho_{\delta}\in W^{1, p}\left(U\right)$. 
Thus, for every $n\in\Nset$ there exists $\phi_{n}\in \Cinfinito\left(U\right)$ such that 
$\abs{\phi_{n} - \ln\rho_{\delta}}_{W^{1,p}} <1/n$, $\phi_{n}\to \ln\rho_{\delta}$ pointwise a.e. and $\ln\delta\leq\phi_n$\footnote{Reminding that the Sobolev space $W^{1,p}\left(\Omega\right)$ is the completion of the set
\be \left\{u\in\Cinfinito\left(\Omega\right) : \int_{\Omega}\abs{u}^{p}dv_{g}<\infty\ \mathrm{and\ } \int_{\Omega}\abs{\nabla u}^{p}dv_{g}<\infty\right\} \nonumber\ee 
with respect to the norm
\be \abs{u}_{W^{1,p}} =\left(\int_{\Omega}\abs{u}^{p}dv_{g} + \int_{\Omega}\abs{\nabla u}^{p}dv_{g}\right)^{1/p}, \nonumber\ee
the approximation result follows by slight modification of classical 
arguments that the reader can find, for instance, in \cite{Au}.}.

Setting $\psi_{n} = e^{\phi_{n}}$ we have that $\psi_{n}\in \Cinfinito\left(U\right)$, $\delta\leq\psi_{n}$, $\psi_n \rightarrow\rho_{\delta}$ a.e. and
\be \int_{K}  \abs{\ln\psi_{n} - \ln\rho_{\delta}}^{p} dv_{g} \longrightarrow 0,  \quad \int_{K} \abs{\frac{\nabla{\psi_{n}}}{\psi_{n}} - 
\frac{\nabla{\rho_{\delta}}}{\rho_{\delta}}}^{p} dv_{g} \longrightarrow 0, 
\quad\mathrm{(as\ } n\rightarrow 0+\infty). \label{limit1}\ee

For every $n\in\Nset$, the function $\varphi_{n}$ defined as 
$\varphi_{n} \decl \frac{\varphi}{\psi_{n}^{p-1}}$ 
belongs to $\Cuno_{0}\left(\Omega\right)$ 
and it is nonnegative since $\varphi\in\Cuno_{0}\left(\Omega\right)$ is 
nonnegative and $\psi_{n}>0$. 
Using $\varphi_{n}$ as test function in $(\ref{subarmonica10})$ we have
\be 0\le \int_{\Omega} \abs{\nabla\rho}^{p-2} \left( \nabla\rho \cdot \nabla\varphi_{n}\right) dv_{g} = \int_{\Omega} \abs{\nabla\rho}^{p-2} \left( \nabla\rho\cdot\nabla\left(\frac{\varphi}{\psi_{n}^{p-1}} \right) \right) dv_{g},  \label{subarmonica12}\ee
which, since by computation  $\nabla\left(\frac{\varphi}{\psi_{n}^{p-1}}\right)=\frac{\nabla\varphi}{\psi_{n}^{p-1}} - \left(p-1\right)\frac{\nabla\psi_{n}}{\psi_{n}^{p}}\varphi$, implies
\be \left(p-1\right) \int_{\Omega}\frac{\abs{\nabla\rho}^{p-2} \nabla\rho\cdot\nabla\psi_{n}}{\psi_{n}^{p}}\, \varphi \ dv_{g}\leq 
\int_{\Omega}\left( \frac{\abs{\nabla\rho}^{p-2}\nabla\rho}{\psi_{n}^{p-1}} \cdot\nabla\varphi \right) dv_{g}.
\label{subarmonica13}\ee
 
Now, letting $n\to +\infty$ we obtain by dominated convergence:
\be \int_{\Omega}\left( \frac{\abs{\nabla\rho}^{p-2}\nabla\rho}{\psi_{n}^{p-1}} \cdot\nabla\varphi\right)  dv_{g} \longrightarrow  \int_{\Omega} \left( \frac{\abs{\nabla\rho}^{p-2}\nabla\rho}{\rho_{\delta}^{p-1}}\cdot\nabla\varphi\right)  dv_{g}, \nonumber\ee
because $\vert\frac{\abs{\nabla\rho}^{p-2}\nabla\rho\cdot\nabla\varphi}{\psi_{n}^{p-1}} \vert
\leq \frac{C\abs{\nabla\rho}^{p-1}}{\delta^{p-1}}\in L^{1}\left(U\right)$. 
Now we claim that 
\be \int_{\Omega} \frac{\abs{\nabla\rho}^{p-2} \nabla\rho\cdot\nabla\psi_{n}}{\psi_{n}^{p}}\, \varphi \, dv_{g} 
=\int_{\Omega} \frac{\abs{\nabla\rho}^{p-2}\nabla\rho}{\psi_{n}^{p-1}}\, \frac{\nabla\psi_n}{\psi_n} \, \varphi \,  dv_{g} 
\longrightarrow  \int_{\Omega} \frac{\abs{\nabla\rho}^{p}}{\rho_{\delta}^{p}}\, \varphi\, dv_{g}. \nonumber\ee
Indeed, 
\be \frac{\abs{\nabla\rho}^{p-2}\nabla\rho}{\psi_{n}^{p-1}} \rightarrow \frac{\abs{\nabla\rho}^{p-2}\nabla\rho}{\rho_{\delta}^{p-1}}     \qquad \mathrm{pointwise\ a.\ e.} \nonumber\ee
and, since $\frac{\abs{\nabla\rho}^{p-2}\nabla\rho}{\psi_{n}^{p-1}} \leq 
\frac{\abs{\nabla\rho}^{p-1}}{\delta^{p-1}}\in L^{p'}\left(U\right)$, 
by Lebesgue dominated convergence theorem we have that 
\be \frac{\abs{\nabla\rho}^{p-2}\nabla\rho}{\psi_{n}^{p-1}} \rightarrow \frac{\abs{\nabla\rho}^{p-2}\nabla\rho}{\rho_{\delta}^{p-1}}     \qquad \mathrm{in\ } L^{p'}\left(U\right). \nonumber\ee
From this and the fact that 
\be  \frac{\nabla\psi_n}{\psi_n}\rightarrow\frac{\nabla\rho}{\rho_{\delta}}  \qquad \mathrm{in\ } L^{p}\left(U\right) \nonumber\ee
we get the claim. 
Therefore, letting $n\rightarrow +\infty$ in $(\ref{subarmonica13})$, we have 
\be \left(p-1\right) \int_{\Omega}\frac{\abs{\nabla\rho}^{p}}{\rho_{\delta}^{p}}\, \varphi \ dv_{g}\leq \int_{\Omega}\left( \frac{\abs{\nabla\rho}^{p-2}\nabla\rho}{\rho_{\delta}^{p-1}} \cdot\nabla\varphi \right) dv_{g},
\nonumber\ee
which is exactly $(\ref{subarmonica11})$, since $\nabla\rho_\delta = \nabla\rho$. 

An application of Lemma $\ref{lemma:k1}$ gives 
\be \left(\frac{p-1}{p}\right)^{p}\int_{\Omega}\frac{\abs{u}^{p}}{\rho_{\delta}^{p}}\abs{\nabla \rho_\delta}^{p} dv_{g}\leq
\int_{\Omega}\abs{\nabla u}^{p} dv_{g}.  \label{limdelta}\ee
Finally, letting $\delta\rightarrow 0$ in $(\ref{limdelta})$ 
and using Fatou's Lemma, we conclude the proof. 
\hfill$\Box$

\section{Further inequalities}\label{sec:gen}

In this section we shall present some slight but natural extensions of Theorem \ref{teo:g}
and  Lemma \ref{lemma:k1}. 
As byproducts of  these generalizations we shall obtain
Hardy inequalities with a weight in the right hand side, Caccioppoli-type inequalities, weighted Gagliardo-Nirenberg inequalities and the uncertain principle.

A first example of a possible generalization of Theorem \ref{teo:g} is the following:

\bt\label{hardygen}
Let $\alpha\in\Rset$, and let $\rho\in W^{1, p}_{loc}\left(\Omega\right)$ be a nonnegative function 
satisfying the following properties:
\begin{enumerate} 
\item[i)]$-\left(p-1-\alpha\right)\Delta_{p}\rho\geq 0$ on $\Omega$ in weak sense,
\item[ii)]$\frac{\abs{\nabla\rho}^{p}}{\rho^{p-\alpha}}, \rho^{\alpha}\in L^{1}_{loc}\left(\Omega\right).$ 
\end{enumerate}
Then the following Hardy inequality holds 
\be \left(\frac{\abs{p-1-\alpha}}{p}\right)^{p}\int_{\Omega}\rho^{\alpha}\frac{\abs{u}^{p}}{\rho^{p}}\abs{\nabla \rho}^{p} dv_{g} \leq \int_{\Omega}\rho^{\alpha}\abs{\nabla u}^{p} dv_{g}, 
\quad u\in\Cinfinito_{0}\left(\Omega\right).
\label{difdis}\ee
\et

The proof of the above theorem is similar to the one of Theorem \ref{teo:g} and it is based on a careful 
choice of the vector field  $h$ and of the function $A_{h}$  in Lemma $\ref{lemma:k1}$.

\noindent{\bf Proof.} 
Let $0<\delta <1$, and $\rho_\delta := \rho + \delta$. 
In order to apply Lemma \ref{lemma:k1} 
we choose the vector field $h$ and the function $A_{h}$ as
\be  h := -\left(p-1-\alpha\right)\frac{\abs{\nabla\rho_\delta}^{p-2}\nabla\rho_\delta}{\rho_{\delta}^{p-1-\alpha}}, \qquad A_{h} := \left(p-1-\alpha\right)^{2}\frac{\abs{\nabla\rho_\delta}^{p}}{\rho_{\delta}^{p-\alpha}}. \label{halfa}\ee
Arguing as in the proof of Theorem $\ref{teo:g}$, 
we have to show that 
\be 
\left(p-1-\alpha\right)^{2} 
\int_{\Omega}\frac{\abs{\nabla\rho_\delta}^{p}}{\rho_{\delta}^{p-\alpha}}\, \varphi \ dv_{g}\leq
\left( p-1-\alpha\right)  \int_{\Omega}  \left( \frac{\abs{\nabla\rho_\delta}^{p-2}\nabla\rho_\delta}{\rho_{\delta}^{p-1-\alpha}}\cdot\nabla\varphi \right) dv_{g},
\label{psub}
\ee 
for every nonnegative function $\varphi\in\C^{1}_{0}\left(\Omega\right)$. 
Let $K\decl supp\varphi \subset\Omega$ and let $U\subset\subset\Omega$ be a neighborhood of $K$. 
Let $k>\delta$, and define $\rho_{k\delta}\decl \inf \left\lbrace \rho_{\delta}, k\right\rbrace$. 
Arguing as in the proof of Theorem $\ref{teo:g}$, 
we have that there exists a sequence $\left\lbrace \psi_n\right\rbrace$ such that $\delta\leq\psi_{n}\leq k$, and 
\be \int_{K}  \abs{\ln\psi_{n} - \ln\rho_{k\delta}}^{p} dv_{g} \longrightarrow 0,  \quad   \int_{K} \abs{\frac{\nabla{\psi_{n}}}{\psi_{n}} - 
\frac{\nabla{\rho_{k\delta}}}{\rho_{k\delta}}}^{p} dv_{g} \longrightarrow 0, 
\quad\mathrm{(as\ } n\rightarrow +\infty). \label{plim}\ee
Then we use $\varphi_{n}:=\frac{\varphi}{\psi_{n}^{p-1-\alpha}}$ as test function in the hypothesis i), 
obtaining 
\be \left(p-1-\alpha\right)^{2} \int_{\Omega}\frac{\abs{\nabla\rho}^{p-2}\nabla\rho\cdot\nabla\psi_{n}}{\psi_{n}^{p-\alpha}}\, \varphi \ dv_{g}\leq
\left( p-1-\alpha\right)  \int_{\Omega}  \left( \frac{\abs{\nabla\rho}^{p-2}\nabla\rho}{\psi_{n}^{p-1-\alpha}}\cdot\nabla\varphi \right) dv_{g}.
\label{psub1} 
\ee 
In the case $\alpha < p-1$ we obtain $(\ref{psub})$ from $(\ref{psub1})$ 
by slight modifications of the proof of Theorem $\ref{teo:g}$, so we will omit the proof. 

Let $\alpha >p-1$. 
We claim that, letting $n\rightarrow +\infty$ in $(\ref{psub1})$, 
and eventually taking a subsequence, we get 
\be\left(p-1-\alpha\right)^{2} \int_{\Omega}\frac{\abs{\nabla\rho}^{p-2}\nabla\rho\cdot\nabla\rho_{k\delta}}{\rho_{k\delta}^{p-\alpha}}\, \varphi \ dv_{g}\leq
\left( p-1-\alpha\right)  \int_{\Omega}  \left( \frac{\abs{\nabla\rho}^{p-2}\nabla\rho}{\rho_{k\delta}^{p-1-\alpha}}\cdot\nabla\varphi \right) dv_{g}.
\label{psub2} 
\ee
In fact, for the right hand side the limit follows by dominated convergence, 
since 
\be \abs {\frac{\abs{\nabla\rho}^{p-2}\nabla\rho}{\psi_{n}^{p-1-\alpha}}\cdot\nabla\varphi} =
\abs{\nabla\rho}^{p-2}\vert\nabla\rho\cdot\nabla\varphi \vert \psi_{n}^{\alpha-p+1}\leq 
C \abs{\nabla\rho}^{p-1}k^{\alpha-p+1}\in L^{1}\left(U\right).  \nonumber\ee
Dealing with the left hand side of $(\ref{psub1})$, we set 
\be \frac{\abs{\nabla\rho}^{p-2}\nabla\rho\cdot\nabla\psi_{n}}{\psi_{n}^{p-\alpha}}  = 
\frac{\abs{\nabla\rho}^{p-2}\nabla\rho}{\psi_{n}^{p-1}} \cdot 
\frac{\nabla\psi_n}{\psi_n} \, \psi_n^{\alpha}  =: f_n \cdot g_n.   \nonumber\ee
As in the proof of Theorem $\ref{teo:g}$, we have 
\be   
f_n \rightarrow \frac{\abs{\nabla\rho}^{p-2}\nabla\rho}{\rho_{k\delta}^{p-1}}   \qquad \mathrm{in\ } L^{p'}\left(U\right), \label{fn}\ee
while from the relations 
\be   \vert g_n\vert \leq \abs{\frac{\nabla\psi_n}{\psi_n}} \cdot k^{\alpha} 
\rightarrow  \abs{\frac{\nabla\rho_{k\delta}}{\rho_{k\delta}}} \cdot k^{\alpha}   \qquad \mathrm{in\ } L^{p}\left(U\right) \nonumber\ee
we obtain that the sequence $g_n$ is bounded in $L^{p}\left(U\right)$. 
Therefore, up to a subsequence, $g_n$ is weakly convergent in $L^{p}\left(U\right)$. 
Since 
\be  g_n \rightarrow  \frac{\nabla\rho_{k\delta}}{\rho_{k\delta}} \cdot \rho_{k\delta}^{\alpha}  
\qquad   \mathrm{pointwise\ a.\ e.}, \nonumber\ee
we have that the convergence is in the weak sense. 
This fact with $(\ref{fn})$ concludes the claim.

Next step is letting $k\rightarrow +\infty$ in $(\ref{psub2})$. 
Let us rewrite the integrand in the right hand side as 
\be   
\abs{ \frac{\abs{\nabla\rho}^{p-2}\nabla\rho}{\rho_{k\delta}^{p-1-\alpha}}\cdot\nabla\varphi } 
= \abs{ \abs{\nabla\rho}^{p-2}\nabla\rho \left(\rho_{k\delta}\right)^{\frac{\alpha -p}{p'}} 
\left(\rho_{k\delta}\right)^{\frac{\alpha -p}{p} +1} \cdot\nabla\varphi  }
\leq  C \abs{\nabla\rho}^{p-1} \left(\rho_{\delta}\right)^{\frac{\alpha -p}{p'}} 
\left(\rho_{\delta}\right)^{\frac{\alpha}{p}},  \nonumber
\ee
which is in $L^{1}\left(U\right)$, since 
$C\abs{\nabla\rho}^{p-1} \left(\rho_{\delta}\right)^{\frac{\alpha -p}{p'}}\in L^{p'}\left(U\right)$ and 
$\left(\rho_{\delta}\right)^{\frac{\alpha}{p}}\in L^{p}\left(U\right)$ by hypothesis ii). 
Thus we can use the dominated convergence to obtain the limit for the right hand side of $(\ref{psub2})$. 

In order to pass to the limit for $k\rightarrow +\infty$ 
in the left hand side of $(\ref{psub2})$, we rewrite the integrand as 
\be   \frac{\abs{\nabla\rho}^{p-2}\nabla\rho\cdot\nabla\rho_{k\delta}}{\rho_{k\delta}^{p-\alpha}}\, \varphi  =   \frac{\abs{\nabla\rho}^{p-2}\nabla\rho\cdot\nabla\rho_{\delta}}{\rho_{k\delta}^{p-\alpha}}\, \chi_{\left\lbrace\rho_{\delta}\leq k\right\rbrace}\varphi   =
\frac{\abs{\nabla\rho}^{p}}{\rho_{k\delta}^{p-\alpha}}  \, \chi_{\left\lbrace\rho_{\delta}\leq k\right\rbrace}\varphi ,   \label{integrand}
\ee
where we have used the fact that $\nabla\rho_{\delta}=\nabla\rho$. 
Now, if $\alpha\leq p$ we apply the dominated convergence, 
since the term in $(\ref{integrand})$ is dominated by the function 
$C\frac{\abs{\nabla\rho}^{p}}{\delta^{p-\alpha}}\in L^{1}\left(U\right)$. 
Whereas, if $\alpha >p$ we use the monotone convergence, 
since $\abs{\nabla\rho}^{p}\rho_{k\delta}^{\alpha -p}\chi_{\left\lbrace\vert\rho_{\delta}\vert\leq k\right\rbrace}\varphi$ is an increasing sequence of nonnegative functions. 

Thus, letting $k\rightarrow +\infty$ in $(\ref{psub2})$, we get 
\be\left(p-1-\alpha\right)^{2} 
\int_{\Omega}\frac{\abs{\nabla\rho}^{p}}{\rho_{\delta}^{p-\alpha}}\, \varphi \ dv_{g}\leq
\left( p-1-\alpha\right)  \int_{\Omega}  \left( \frac{\abs{\nabla\rho}^{p-2}\nabla\rho}{\rho_{\delta}^{p-1-\alpha}}\cdot\nabla\varphi \right) dv_{g}, 
\label{psub3} 
\ee
which is exactly $(\ref{psub})$, since $\nabla\rho_\delta = \nabla\rho$. 

As in Theorem \ref{teo:g}, an application of lemma \ref{lemma:k1} gives 
\be \left(\frac{\alpha -p+1}{p}\right)^{p}\int_{\Omega}\rho_{\delta}^{\alpha}\frac{\abs{u}^{p}}{\rho_{\delta}^{p}}\abs{\nabla \rho_\delta}^{p} dv_{g} \leq \int_{\Omega}\rho_{\delta}^{\alpha}\abs{\nabla u}^{p} dv_{g}.
\label{limdelta1}\ee
Finally, letting $\delta\rightarrow 0$ in $(\ref{limdelta1})$, we conclude the proof. 
Indeed we can use the dominated convergence for the right hand side, since 
$\rho_{\delta}^{\alpha}\abs{\nabla u}^{p}\leq C\left(\rho +\delta\right)^\alpha \leq C\left(\rho +1\right)^\alpha 
\in L^{1}_{loc}\left(U\right)$, and apply the Fatou's Lemma for the left hand side. 
\hfill$\Box$

\boss\label{rem:amp} 
If $\alpha\le p$, then the hypothesis $ii)$ in the above Theorem  $\ref{hardygen}$, 
can be avoided. Indeed, since $\rho\in W^{1,p}_{loc}\left(\Omega\right)$ we get that $\rho^\alpha\in L^1_{loc}\left(\Omega\right)$ and from the proof it follows that $\frac{\abs{\nabla\rho}^{p}}{\rho^{p-\alpha}}\in L^1_{loc}\left(\Omega\right)$.
\eoss

As a consequence of Theorem $\ref{hardygen}$ (it suffices to take $\alpha=p+q$), 
we obtain the following Caccioppoli-type inequality for $p$-subharmonic functions, which is worth of mention:

\bc($L^p$-Caccioppoli-type inequality) \label{cacc}
Let $\rho\in L^{1}_{loc}\left(M\right)$ and $q>-1$. 
Assume that $\rho$ is nonnegative on an open set $\Omega\subset M$, 
and $\rho\in W^{1, p}_{loc}\left(\Omega\right)$, 
$\rho^q \abs{\nabla\rho}^{p}, \rho^{p+q}\in L^{1}_{loc}\left(\Omega\right)$. 
If $\Delta_{p}\rho\geq 0$ on $\Omega$ in weak sense, 
then we have  
\be \left(\frac{q+1}{p}\right)^{p}\int_{\Omega}\rho^q\abs{\nabla \rho}^{p}\abs{u}^{p} dv_{g} \leq\int_{\Omega}\rho^{p+q}\abs{\nabla u}^{p} dv_{g}, 
\quad u\in\Cinfinito_{0}\left(\Omega\right).
\label{ww}\ee
\ec

Notice that for $q=0$ and $p=2$ the above theorem is a version of the 
classical Caccioppoli inequality on manifolds. 
See also \cite{PRS} for a version of Caccioppoli inequality related to subharmonic functions on manifolds.

\bigskip 

Now we present a possible generalization of Lemma \ref{lemma:k1}, and some of its consequences, 
like the weighted Gagliardo-Nirenberg inequality and the uncertain principle on manifolds.

\bl\label{lemma:k2} 
Let $h\in L^{1}_{loc}\left(\Omega\right)$ be a vector field and let 
$A_{h}\in L^{1}_{loc}\left(\Omega\right)$ be a nonnegative function such that
\begin{enumerate} 
\item[i)] $A_{h} \leq\Div h,$
\item[ii)] $\frac{\abs{h}^{p}}{A_{h}^{p-1}} \in L^{1}_{loc}\left(\Omega\right)$. 
\end{enumerate}

Then for every $u\in\Cuno_{0}\left(\Omega\right)$, $q\in \Rset$, $s>0$ and $a>1$  we have
\be\int_{\Omega} \abs{u}^{s} \abs{h}^q dv_{g}\leq p^{p/a} 
\left(\int_{\Omega} \frac{\abs{h}^{p}}{A_{h}^{p-1}} \abs{\nabla{u}}^{p} dv_{g}\right)^{1/a} 
\left(\int_{\Omega} \frac{\abs{h}^{qa'}}{A_{h}^{a'-1}} \abs{{u}}^{\frac{as-p}{a-1}} dv_{g}   \right)^{1/a'}, \label{lem102}\ee
provided $\abs h ^q\in L^1_{loc}(\Omega)$.

In particular, setting $w:=\abs{h}A_{h}^{\frac{1-p}{p}}$, we have
\begin{enumerate}
\item 
\be \left( \int_{\Omega} \abs{u}^{s} \abs{h}^q dv_{g}\right)^{1/s} \leq p^{q(p-1)/s} 
\left(\int_{\Omega} w^p \abs{\nabla{u}}^{p} dv_{g}\right)^{b/p} 
\left(\int_{\Omega} w^{t\delta}\abs{{u}}^{\delta} dv_{g}   \right)^{(1-b)/\delta}, \label{lemgn}\ee
where $t,\delta >0$ and 
$$ \frac{1}{s} = \frac{b}{p} + \frac{1-b}{\delta}, \qquad 
\frac{1}{q} = \frac{1}{p'} + \frac{1}{t\delta},   \qquad 
b = \frac{t(p-1)}{1+t(p-1)}.  $$
\item 
\be\int_{\Omega} \abs{u}^{s}  dv_{g}\leq p^{p/a} 
\left(\int_{\Omega} w^p \abs{\nabla{u}}^{p} dv_{g}\right)^{1/a} 
\left(\int_{\Omega} \frac{1}{A_{h}^{a'-1}} \abs{{u}}^{\frac{as-p}{a-1}} dv_{g}   \right)^{1/a'}, \label{lemph}\ee
where $s>0$ and $a>1$.
\end{enumerate}

\el
\bp By H\"older inequality with exponent $a$ we have
\bern
\int_{\Omega} \abs{u}^{s} \abs{h}^q dv_{g}&=& 
\int_{\Omega} \abs{u}^{p/a} A_h^{1/a}  \abs{h}^q  A_h^{-1/a} \abs{u}^{s-p/a} dv_{g}\\
&\leq & \left(\int_{\Omega} \abs{u}^{p} A_h \ dv_{g} \right)^{1/a}    
\left(\int_{\Omega} \abs{h}^{q a'}  A_h^{-a'/a} \abs{u}^{\frac{as -p}{a-1}} dv_{g} \right)^{1/a'} ,
\eern
which by using (\ref{lem10}), implies (\ref{lem102}).

From (\ref{lem102}) we get (\ref{lemgn}) by choosing $a=1 + \frac{p'}{t\delta}$, 
and (\ref{lemph}) by choosing $q=0$.  
\ep

Specializing $h$ and $A_h$ we obtain from (\ref{lemgn}) and (\ref{lemph})
a weighted Gagliardo-Nirenberg inequality 
and an uncertain principle respectively. 
In particular, choosing $h$ and $A_h$ as in $(\ref{hAh})$, we have the following

\begin{theorem} 
Let $\rho\in W^{1, p}_{loc}\left(\Omega\right)$ be nonnegative.
Assume that $\rho$ is $p$-superharmonic function on $\Omega\subset M$ and satisfies the hypotheses of Theorem \ref{teo:g}. 
Let $\delta >0$ and $0\leq b \leq 1$. 
Then for every $u\in\Cuno_{0}\left(\Omega\right)$, we have
\be\left( \int_{\Omega} \abs{u}^{s} \frac{\abs{\nabla\rho}^{q(p-1)}}{\rho^{q(p-1)}}\, dv_{g}\right)^{1/s} \leq\left(\frac{p}{p-1}\right)^{q(p-1)/s} 
\left(\int_{\Omega}  \abs{\nabla{u}}^{p} dv_{g}\right)^{b/p} 
\left(\int_{\Omega} \abs{{u}}^{\delta} dv_{g}   \right)^{(1-b)/\delta},  \label{lemgnS}\ee
where
$$ 
\frac{1}{s} = \frac{b}{p} + \frac{1-b}{\delta}, \qquad 
\frac{1}{q(p-1)} = \frac{1}{p} + \frac{1-b}{b\delta}.  $$
In particular, if $\rho = d^\alpha$ for some $\alpha \neq 0$ with 
$\abs{\nabla d}=1$, then we have
\be\int_{\Omega} \frac{\abs{u}^{s}}{d^{p-1}}\, dv_{g}\leq
\left(\frac{p}{\abs \alpha (p-1)}\right)^{(p-1)} 
\left(\int_{\Omega}  \abs{\nabla{u}}^{p} dv_{g}\right)^{1/p'} 
\left(\int_{\Omega} \abs{{u}}^{\delta} dv_{g}   \right)^{1/p},  \label{lemgnSS}\ee
where $s=p-1 +\frac{\delta}{p}$.
\end{theorem}

Notice that for $s=p=2$ the inequality (\ref{lemgnSS}) is 
the weighted Gagliardo-Nirenberg inequality on manifold. 
Its counterpart in Euclidean setting is largely studied by many authors, see for instance \cite{DV}.
Further examples of manifolds and functions $\rho$ satisfying the hypotheses of the above 
theorem are given in Section \ref{sec:ex}.

\begin{theorem} 
Let $\rho\in W^{1, p}_{loc}\left(\Omega\right)$ be nonnegative.
Assume that $\rho$ is $p$-superharmonic function on $\Omega\subset M$ and satisfies the hypotheses of Theorem \ref{teo:g}. 
Let $s>0$ and $a>1$. Then for every $u\in\Cuno_{0}\left(\Omega\right)$, we have
\be \int_{\Omega} \abs{u}^{s} dv_{g}\leq\left(\frac{p}{p-1}\right)^{p/a} 
\left(\int_{\Omega}  \abs{\nabla{u}}^{p} dv_{g}\right)^{1/a} 
\left(\int_{\Omega} \abs{{u}}^{\frac{as-p}{a-1} } \frac{\rho^{p(a'-1)}}{\abs{\nabla \rho}^{p(a'-1)}}\, dv_{g}   \right)^{1/a'}. \label{lemphS}\ee

In particular, if $\rho = d^\alpha$ for some $\alpha \neq 0$ with 
$\abs{\nabla d}=1$, then we have \be\int_{\Omega} \abs{u}^{s} dv_{g}\leq 
\left(\frac{p}{\abs \alpha (p-1)}\right)^{p/a} 
\left(\int_{\Omega}  \abs{\nabla{u}}^{p} dv_{g}\right)^{1/a} 
\left(\int_{\Omega} \abs{{u}}^{\frac{as-p}{a-1} } d^{p(a'-1)}dv_{g}   \right)^{1/a'}. \label{lemphSS}\ee
\end{theorem}

Notice that if  $a=s=p=2$ the inequality (\ref{lemphSS}) in the Euclidean setting 
coincides with the celebrated uncertain principle with $d=\abs x$, the Euclidean norm.

\boss 
Different choices of the vector field $h$ and of the function $A_h$ 
in Lemma $\ref{lemma:k2}$, produce inequalities different than (\ref{lemgnS})---(\ref{lemphSS}). 
For instance, one can define $h$ and $A_h$ as in $(\ref{halfa})$,
obtaining a version of   (\ref{lemgnS})---(\ref{lemphSS}) with further weights.
\eoss

\bigskip 

To end this section, we want to point out that it is possible to extend all the results of this paper 
considering vector fields of the type 
$\nabla_{\mu}u := \mu(\nabla u)$, 
where $\mu$ is a $\left(1, 1\right)$-tensor (say  $\Cuno$). 
In this case, replacing $\nabla$ with $\nabla_\mu$, 
a Hardy-type inequality like (\ref{dishardy10}) holds provided
$\nabla^{\ast}_{\mu}\left(\abs{\nabla_{\mu} u}^{p-2}\nabla_{\mu} u\right)\geq 0$, 
where $\nabla^{\ast}_{\mu}$ stands for the adjoint of $\nabla_{\mu}$. 
We leave the details to the interested reader. 
Notice that the study of Hardy inequalities for the vector field $\nabla_\mu$
was already studied in \cite{Da4}, 
when  the support of the manifold is $\Rset^N$.

\section{Remarks on the best constant}\label{sec:bc}

Theorems \ref{teo:g} and \ref{hardygen} affirm the validity of some 
Hardy inequalities with an explicit value of the constants involved.
In many cases these constants, $\left(\frac{p-1}{p}\right)^{p}$ and 
$\left(\frac{\abs{p-1-\alpha}}{p}\right)^{p}$, result to be sharp. 
For example in \cite{Da4} the author proves the sharpness of the constant 
$\left(\frac{\abs{p-1-\alpha}}{p}\right)^{p}$ involved in the 
inequality of Theorem \ref{hardygen} in several cases. 
Moreover the question of the existence of functions 
that realize the best constant is analyzed in many papers 
(for the Euclidean case see for instance \cite{BM, DH, MMP, MS1, MS2}).
On the other hand, the knowledge of the best constants for the inequalities 
plays a crucial role in \cite{AX, CX, Xi}.

For the sake of simplicity, we shall focus our attention on the inequality $(\ref{dishardy10})$. 
We denote by $c\left(\Omega\right)$ the best constant in $(\ref{dishardy10})$, namely
\be c\left(\Omega\right) := \inf_{u\in D^{1,p}\left(\Omega\right), u\neq0} \frac{\int_{\Omega} \abs{\nabla u}^{p} dv_{g}}{\int_{\Omega}\frac{\abs{u}^{p}}{\rho^{p}}\abs{\nabla \rho}^{p} dv_{g}}.\label{best10}\ee
Then, we have the following:

\bt\label{teo:g1} Under the same hypotheses of Theorem \ref{teo:g} we have: 
\begin{enumerate}
\item[1.] If $\rho^{\frac{p-1}{p}}\in D^{1,p}\left(\Omega\right)$, 
then $c\left(\Omega\right)=\left(\frac{p-1}{p}\right)^{p}$ and $\rho^{\frac{p-1}{p}}$ is a minimizer.
\item[2.] If $\rho^{\frac{p-1}{p}}\notin D^{1,p}\left(\Omega\right)$, $p\geq2$ and $c\left(\Omega\right)=\left(\frac{p-1}{p}\right)^{p}$, 
then the best constant $c\left(\Omega\right)$ is not achieved.
\end{enumerate}
\et

\noindent{\bf Proof.}
1) From $(\ref{dishardy10})$ we have $c\left(\Omega\right)\geq\left(\frac{p-1}{p}\right)^{p}$. 
Moreover, if $\rho^{\frac{p-1}{p}}\in D^{1,p}\left(\Omega\right)$, by computation
\begin{eqnarray} 
\int_{\Omega} \abs{\nabla\rho^{\frac{p-1}{p}}}^{p} dv_{g} &=& 
\int_{\Omega}\left(\frac{p-1}{p}\right)^{p} \left(\rho^{-1/p}\right)^{p} \abs{\nabla\rho}^{p} dv_{g} \nonumber\\
&=& \left(\frac{p-1}{p}\right)^{p} \int_{\Omega} \frac{\abs{\nabla\rho}^{p}}{\rho}\, dv_{g} \nonumber\\
&=&  \left(\frac{p-1}{p}\right)^{p} \int_{\Omega} \frac{\abs{\rho^{\frac{p-1}{p}}}^{p}}{\rho^{p}}\abs{\nabla \rho}^{p} dv_{g}.    \nonumber\end{eqnarray}
Thus, taking $u = \rho^{\frac{p-1}{p}}$, we obtain the infimum in $(\ref{best10})$.

2) Let $u\in\Cinfinito_{0}\left(\Omega\right)$. We define the functional $I$ as
\be I\left(u\right) := \int_{\Omega} \abs{\nabla u}^{p} dv_{g} - \left(\frac{p-1}{p}\right)^{p} \int_{\Omega} \frac{\abs{u}^{p}}{\rho^{p}}\abs{\nabla \rho}^{p} dv_{g}. \nonumber \ee
We note that the functional $I$ is nonnegative, 
since $(\ref{dishardy10})$ holds, and the best constant will be achieved if and only if 
$I\left(u\right)=0$ for some $u\in D^{1, p}\left(\Omega\right)$.

Let $v$ be the new variable $v := \frac{u}{\rho^{\gamma}}$ with $\gamma := \frac{p-1}{p}$. 
By computation we have
\begin{eqnarray}
\abs{\nabla u}^{2} &=&  
\abs{\nabla \left(v\rho^{\gamma}\right)}^{2} \nonumber\\
&=& \abs{\gamma}^{2} v^{2} \rho^{2\gamma-2}\abs{\nabla\rho}^{2} 
+ \rho^{2\gamma} \abs{\nabla v}^{2} 
+ 2\gamma v \rho^{2\gamma-1} \nabla\rho\cdot\nabla v. \label{inequ}
\end{eqnarray}
(If $\rho$ is not smooth enough, we can consider $\psi_{n}$ as in the proof 
of Theorem $\ref{teo:g}$ and after the computation take the limit as $n\rightarrow +\infty$).

We remind that the inequality
\be \left(\xi-\eta\right)^{s} \geq \xi^{s} - s\eta\xi^{s-1} \label{ineq}\ee
holds for every $\xi, \eta, s\in\Rset$, with $\xi>0$, $\xi>\eta$ and $s\geq1$ (see \cite{GGM}). 
Applying $(\ref{ineq})$ and $(\ref{inequ})$, with $s = p/2$, $\xi = \abs{\gamma}^{2} v^{2} \rho^{2\gamma-2}\abs{\nabla\rho}^{2}$ and $\eta = - 2\gamma v \rho^{2\gamma-1} \nabla\rho\cdot\nabla v - \rho^{2\gamma} \abs{\nabla v}^{2}$, 
we have 
\bern 
\abs{\nabla u}^{p}\geq \abs{\gamma}^{p} v^{p} \rho^{p\gamma-p} \abs{\nabla\rho}^{p} 
+ p\abs{\gamma}^{p-2}\gamma \abs{v}^{p-2} v \abs{\nabla\rho}^{p-2}\nabla\rho\cdot\nabla v\\
+ \frac{p}{2}\abs{\gamma}^{p-2}  \abs{v}^{p-2} \rho \nabla\rho^{p-2}  \abs{\nabla v}^{2}.
\eern
Then, taking into account that $v = \frac{u}{\rho^{\gamma}}$, we have
\begin{eqnarray}
I(u)&=& \int_{\Omega}\abs{\nabla u}^{p} dv_{g} - \abs{\gamma}^{p} \int_{\Omega}\frac{\abs{u}^{p}}{\rho^{p}}\abs{\nabla\rho}^{p} dv_{g}\\ &\geq& \int_{\Omega} p\abs{\gamma}^{p-2}\gamma \abs{v}^{p-2} v \abs{\nabla\rho}^{p-2} \left( \nabla\rho\cdot\nabla v\right)  dv_{g} \nonumber\\
&&\ \ \ + \int_{\Omega}\frac{p}{2}\abs{\gamma}^{p-2}  \abs{v}^{p-2} \rho \abs{\nabla\rho}^{p-2}  \abs{\nabla v}^{2}  dv_{g}=: I_{1}\left(v\right)  + I_{2}\left(v\right). \label{eps10}\end{eqnarray}

Re-arranging the expression in $I_{1}\left(v\right)$ and integrating by parts we obtain
\bern
I_{1}\left(v\right)  &=& \left(\frac{p-1}{p}\right)^{p-1} \int_{\Omega} \nabla\left(\abs{v}^{p}\right) 
\cdot\abs{\nabla\rho}^{p-2}\nabla\rho \  dv_{g} \\
&=&  \left(\frac{p-1}{p}\right)^{p-1}  \int_{\Omega} \abs{v}^{p} \left(-\Delta_{p}\rho\right) dv_{g} \geq 0, \nonumber\eern
where we have used the hypothesis $-\Delta_{p}\rho\geq 0$. 
On the other hand we can rewrite $I_{2}\left(v\right)$ as
\be I_{2}\left(v\right) =  \frac{2}{p}\abs{\gamma}^{p-2} \int_{\Omega} \rho \abs{\nabla\rho}^{p-2} \abs{\nabla\abs{v}^{p/2}}^{2} dv_{g}.\ee

\noindent Thus, we conclude that for every $u\in D^{1, p}\left(\Omega\right)$
\be I\left(u\right)  \geq \frac{2}{p}\abs{\gamma}^{p-2} \int_{\Omega} \rho \abs{\nabla\rho}^{p-2} \abs{\nabla\abs{v}^{p/2}}^{2} dv_{g}>0 , \nonumber\ee
and this inequality implies the non existence of minimizers in $D^{1,p}\left(\Omega\right)$.
\hfill$\Box$

\bigskip
We end this section by showing a further result that arises from the fact 
that the best constant, in some cases, is not achieved. 
Indeed, if the best constant involved in an inequality is not achieved, 
it is natural to ask if a reminder term can be added. 
The next result shows that in the inequality $(\ref{dishardy10})$ 
one can add a reminder term.
\begin{theorem}\label{teo:resto} Let $p=2$ and let $\rho$ be as Theorem \ref{teo:g}. We  define
$$ \Lambda_1\decl \inf_{u\in C_0^1(\Omega)} \frac{\int_\Omega \rho \abs{\nabla u}^2 dv_g}{\int_\Omega \rho \abs{ u}^2 dv_g}.$$
Assume that $\Lambda_1>0$. Then
\be  \int_{\Omega} \abs{\nabla u}^{2} dv_{g} \ge \frac{1}{4} 
\int_{\Omega} \frac{\abs{u}^{2}}{\rho^{2}}\abs{\nabla \rho}^{2} dv_{g} 
+\Lambda_1 \int_{\Omega} u^2 \,dv_{g}, \qquad u\in\Cuno_{0}\left(\Omega\right). \label{hresto} \ee
\end{theorem}
\bp 
We shall give a  sketch of the proof since it is  similar to the proof 
of Theorem \ref{teo:g1}. By using the same notation of the proof of Theorem \ref{teo:g1}, 
from (\ref{inequ}), we deduce that
$$ I(u)\ge \int_\Omega \rho \abs{\nabla v}^2\,dv_{g}\ge 
\Lambda_1\int_\Omega \rho \abs{ v}^2 dv_g=\Lambda_1\int_\Omega  u^2\,dv_g,$$
where we have used the fact that $-\Delta \rho\ge 0$, 
the definition of $\Lambda_1$ and  $v=u/\rho^{1/2}$. This concludes the proof.
\ep

An example of manifold where Theorem \ref{teo:resto} applies is 
the following. Let $\rho$ be a nonnegative superharmonic function on $\Rset^N$
and let $\Omega\subset \Rset^N$ a bounded open set. Then $\rho$ belongs to
the Muckenhoupt class $A_1$ and this implies that $\Lambda_1>0$ 
(indeed it suffices to combine Theorems 3.59 and 15.21 of \cite{HKM}). 
In particular, with the choice $\rho\decl \abs x^{2-N}$, $N>2$ and $\Omega\subset \Rset^N$ a bounded open set, $(\ref{hresto})$ reads as 
$$  \int_{\Omega} \abs{\nabla u}^{2} dx\ge \frac{(N-2)^2}{4} 
\int_{\Omega} \frac{\abs{u}^{2}}{\vert x\vert^{2}}\, dx 
+\Lambda_1 \int_{\Omega} u^2 \,dx, \qquad u\in\Cuno_{0}\left(\Omega\right), $$
which is the celebrated inequality proved in \cite{BV}. 
See also \cite{Da3, GGM} for related results in Euclidean and subelliptic setting for $p>1$ and for further references.

\section{First order interpolation inequalities}\label{sec:ckn}

In this section we shall study some inequalities of Hardy-Sobolev type. 
As already said above, interpolation inequalities
as well as the knowledge of an estimate of the best constant 
have an important role in several areas of mathematical science. 
Thus we shall address some efforts to keep track of explicit values 
of the involved constants. 

We shall assume that the Sobolev inequality 
\be S\left(p\right)\left(\int_{\Omega}\abs{u}^{p^{*}}dv_{g}\right)^{1/p^{*}} 
\leq\left(\int_{\Omega}\abs{\nabla u}^{p} dv_{g}\right)^{1/p}, \quad u\in\Cinfinito_{0}\left(\Omega\right) \eqno{(S)}\nonumber\ee
holds for some $p^{*}>0$, 
and the Hardy inequality 
\be H\left(\alpha, p\right) \int_{\Omega}\rho^{\alpha}\frac{\abs{u}^{p}}{\rho^{p}}\abs{\nabla \rho}^{p} dv_{g} \leq \int_{\Omega}\rho^{\alpha}\abs{\nabla u}^{p} dv_{g}, \quad u\in\Cinfinito_{0}\left(\Omega\right)  
\eqno{(H_{\alpha})}\nonumber\ee
holds for an exponent $\alpha\in\Rset$. 

In some cases, the validity of $(S)$ implies that $(H_\alpha)$ holds as well.
Indeed, let $N>2$ and let $M$ be a $N$-dimensional complete and connected Riemannian manifold  with infinite volume, if $(S)$ holds with $p=2$ and $p^*=2N/(N-2)$ then
$M$ is hyperbolic (see \cite{Ca2}). 
In this case, from Theorem \ref{teo:parabolic}
we have that a Hardy inequality holds. Therefore, there exists a nonnegative
non constant superharmonic function $\rho$ ad hence $(H_\alpha)$ holds with
$p=2$ and $\alpha<1$ (see Theorem \ref{hardygen} and Remark \ref{rem:amp}).

\medskip 

In order to state our main result of this section, we need the following preliminary theorem:

\bt\label{teo:hardy-sob} 
Assume that $(S)$ holds on $\Omega$. 
Let $\theta\in\Rset$ and $\rho\geq0$ be a function such that $(H_{\alpha})$ holds with $\alpha=p\theta$. 
Then there exists $C_2>0$ such that 
\be C_2 \left(\int_{\Omega} \rho^{p^{*}\theta} \abs{u}^{p^{*}}dv_{g}\right)^{1/p^{*}}\leq 
\left(\int_{\Omega}\rho^{p\theta}\abs{\nabla u}^{p} dv_{g}\right)^{1/p}, \quad u\in\Cinfinito_{0}\left(\Omega\right). \label{hsob}\ee 
Moreover 
\be    C_2= 
S\left(p\right)\frac{H\left(p\theta,p\right)^{1/p}}
{\abs \theta + H\left(p\theta,p\right)^{1/p}}. 
\nonumber\ee 

In particular, if $\rho\in L^{1}_{loc}\left(\Omega\right)$ is 
a nonnegative function satisfying the hypotheses of Theorem $\ref{hardygen}$ with $\alpha=p\theta$ 
and $(S)$ holds, then we obtain $(\ref{hsob})$ with 
\be   C_2= S(p)
\frac{\abs{p-1-p\theta}}{p\abs\theta +\abs{p-1-p\theta}}
\nonumber\ee
\et

\noindent{\bf Proof.} 
The case $\theta=0$ corresponds to the Sobolev inequality. Let $\theta\neq 0$. 
Let $u\in\Cinfinito_{0}\left(\Omega\right)$ and define $v$ as $v := \rho^{\theta}u$. 
By computation we have 
\begin{eqnarray}
\abs{\nabla v}^{p} &=&  \abs{\nabla \left(\rho^{\theta}u\right)}^{p} 
= \abs{\rho^{\theta}\nabla u + \theta\rho^{\theta-1}u\nabla\rho}^{p} \nonumber\\
&\leq& \left(\rho^{\theta}\abs{\nabla u} + \abs{\theta}\rho^{\theta-1}\abs{u}\abs{\nabla\rho}\right)^{p} \nonumber\\ 
&=&\left(\rho^{\theta}\abs{\nabla u} + \frac{\abs{\theta}}{H^{1/p}}H^{1/p}\rho^{\theta-1}\abs{u}\abs{\nabla\rho}\right)^{p} \label{cal}
\end{eqnarray}
where, for sake of brevity, $H= H\left(p\theta,p\right)$ and $S=S(p)$. 
By using the inequality
$$ (a+\frac{1-\epsilon}{\epsilon}b)^p\le \epsilon^{1-p} a^p + 
\frac{1-\epsilon}{\epsilon^p}b^p\qquad (0<\epsilon<1,\ \  a,b>0)$$
with $\epsilon\decl \frac{H^{1/p}}{H^{1/p}+\abs \theta}$,
$a\decl  \rho^{\theta}\abs{\nabla u}$ and 
$b\decl H^{1/p}\rho^{\theta-1}\abs{u}\abs{\nabla\rho}$, we have
\be  \abs{\nabla v}^{p}\leq \epsilon^{1-p}
\rho^{p\theta}\abs{\nabla u}^{p} + 
\frac{1-\epsilon}{\epsilon^p} H \rho^{p\theta-p}\abs{u}^{p}\abs{\nabla\rho}^{p}.   \label{cal2} \ee

Then, by $(S)$ and using $(H_{\alpha})$ with $\alpha=p\theta$, we obtain
\begin{eqnarray}\left(\int_{\Omega}\rho^{p^* \theta}\abs{u}^{p^{*}}dv_{g}\right)^{p/p^{*}}
&=&\left(\int_{\Omega}\abs{v}^{p^{*}}dv_{g}\right)^{p/p^{*}} \le 
S^{-p} \int_{\Omega}\abs{\nabla v}^{p}  \\
&\leq& S^{-p}  \left[ \epsilon^{1-p} \int_{\Omega}\rho^{p\theta}\abs{\nabla u}^{p}dv_{g} +  \frac{1-\epsilon}{\epsilon^p} H  \int_{\Omega} \rho^{p\theta}\frac{\abs{u}^{p}}{\rho^{p}}\abs{\nabla\rho}^{p} dv_{g}\right] \nonumber\\
&\leq& S^{-p}\left[\epsilon^{1-p} + \frac{1-\epsilon}{\epsilon^p}\right] 
\int_{\Omega}\rho^{p\theta}\abs{\nabla u}^{p}dv_{g} = (S\epsilon)^{-p}\int_{\Omega}\rho^{p\theta}\abs{\nabla u}^{p}dv_{g},\nonumber
\end{eqnarray}
which concludes the proof.
\hfill$\Box$

\bigskip

\bt \label{cknmanif}
Assume that $(S)$ holds on $\Omega$ with $p^*>p$.
Let $\theta\in\Rset$ and $\rho\geq0$ be a function such that $(H_{\alpha})$ holds with $\alpha=p\theta$. 
Let $r>0$, $0\leq a\leq1$, $\gamma$, $\epsilon, \sigma$ and $\delta$ be real numbers 
satisfying the following relations
\be \frac{1}{p}\ge\frac 1r\ge\frac{1-a}{p}+\frac{a}{p^*},  \label{condr}\ee
\be \gamma + \frac{p^{*}\left(r-p\right)}{r\left(p^{*}-p\right)}=(1-\theta) a+\delta\left(1-a\right)  \label{cond1}\ee
and
\be \epsilon=\theta a +\sigma\left(1-a\right). \label{cond2}\ee
Then there exists $C_3>0$ such that  
\begin{eqnarray} 
C_3\left(\int_{\Omega}\frac{\abs{u}^{r}}{\rho^{\gamma r}}\abs{\nabla \rho}^{\left(\gamma+\epsilon\right)r} dv_{g}\right)^{1/r} 
&\leq&  \left(\int_{\Omega}\rho^{\theta p}\abs{\nabla u}^{p} dv_{g}\right)^{a/p}  \left(\int_{\Omega}\abs{u}^{p}\frac{\abs{\nabla\rho}^{(\delta+\sigma)p}}{\rho^{\delta p}}\, dv_{g}\right)^{\left(1-a\right)/p} , \nonumber\\
&&\qquad\qquad\qquad\qquad\qquad u\in\Cinfinito_{0}\left(\Omega\right), 
\label{disgen}\end{eqnarray}
that is 
\be C_3\left|u \frac{\abs{\nabla \rho}^{\gamma+\epsilon}}{\rho^{\gamma}}\right|_{L^{r}} \leq  \left|\rho^{\theta}\abs{\nabla u}\right|^{a}_{L^{p}}\left|u\frac{\abs{\nabla\rho}^{\delta+\sigma}}{\rho^{\delta}}\right|^{1-a}_{L^{p}}, \ee
\noindent provided $\frac{\abs{\nabla\rho}^{p\sigma}}{\rho^{p\delta}}\in L^{1}_{loc}\left(\Omega\right)$. 

\noindent Moreover 
\be   C_3 = C_2^{\frac{p^{*}\left(r-p\right)}{r\left(p^{*}-p\right)}}  H\left(p\theta,p\right)^{\frac{a}{p} - \frac{p^{*}\left(r-p\right)}{pr\left(p^{*}-p\right)}}. \nonumber\ee

In particular, if $\rho\in L^{1}_{loc}\left(\Omega\right)$ is 
a nonnegative function satisfying the hypotheses of Theorem $\ref{hardygen}$ with $\alpha=p\theta$ 
and $(S)$ holds, then we obtain $(\ref{disgen})$ with 
\be   C_3 = C_2^{\frac{p^{*}\left(r-p\right)}{r\left(p^{*}-p\right)}}  \left(\frac{\abs{p-1-p\theta}}{p}\right)^{a - \frac{p^{*}\left(r-p\right)}{r\left(p^{*}-p\right)}}. \nonumber\ee
\et

\noindent{\bf Proof.} From condition (\ref{condr}) it follows that 
$p^*\ge r \ge p$. We shall distinguish tree cases.

Case: $r=p^*$. From (\ref{condr}) necessarily we have $a=1$ and hence from 
(\ref{cond1}) and (\ref{cond2}) $\epsilon=\theta=-\gamma$. The inequality to prove is actually the thesis of Theorem \ref{teo:hardy-sob}.

Case: $r=p$. If $a=0$ there is nothing to prove. 
If $a=1$, then the thesis is the the inequality $(H_\alpha)$. 
Let $0<a<1$. By using  (\ref{cond1}) and (\ref{cond2}) we have
$$\int_{\Omega}\frac{\abs{u}^{r}}{\rho^{\gamma r}}\abs{\nabla \rho}^{\left(\gamma +\epsilon\right)r} dv_{g} = 
\int_{\Omega} \frac{\abs{u}^{ap}}{\rho^{(1-\theta) a p}}\abs{\nabla \rho}^{ap}
\frac{\abs{u}^{(1-a)p}}{\rho^{(1-a) \delta p}}
\abs{\nabla \rho}^{\left(\delta +\sigma\right)p} dv_{g}. $$
Now the claim follows applying H\"older inequality with exponent $1/a$ and then Hardy inequality  $(H_\alpha)$.

Case: $p^*>r>p$.
Let $q\in\Rset$ be a parameter that we shall fix later. 
Using H\"older inequality with exponent $s>1$ we obtain
\begin{eqnarray}
\int_{\Omega}\frac{\abs{u}^{r}}{\rho^{\gamma r}}\abs{\nabla \rho}^{\left(\gamma +\epsilon\right)r} dv_{g} =
\int_{\Omega}\abs{u}^{r-q}\rho^{p^{*}\theta/s} \frac{\abs{u}^{q}}{\rho^{\gamma r+p^{*}\theta/s}} \abs{\nabla \rho}^{\left(\gamma + \epsilon\right) r} dv_{g} \nonumber\\
\leq \left(\int_{\Omega}\abs{u}^{\left(r-q\right)s}\rho^{p^{*}\theta} dv_{g}\right)^{1/s} 
\left(\int_{\Omega}\frac{\abs{u}^{qs'}}{\rho^{\left(\gamma r+p^{*}\theta/s\right)s'}}\abs{\nabla \rho}^{\left(\gamma +\epsilon\right)rs'} dv_{g}\right)^{1/s'}. \label{gen1}
\end{eqnarray}
Now we apply H\"older inequality with exponent $t>1$ to the second term of $(\ref{gen1})$ and obtain
\begin{eqnarray} 
\int_{\Omega}\frac{\abs{u}^{qs'}}{\rho^{\left(\gamma r+p^{*}\theta/s\right)s'}}\abs{\nabla \rho}^{\left(\gamma +\epsilon\right)rs'} dv_{g} 
= \int_{\Omega}\frac{\abs{u}^{qs'/t}}{\rho^{qs'/t}} \abs{\nabla \rho}^{qs'/t}\rho^{p\theta/t} 
\frac{\abs{u}^{qs'/t'}\abs{\nabla \rho}^{\left(\gamma + \epsilon\right)rs' - qs'/t}}{\rho^{\left(\gamma r+p^{*}\theta/s\right)s' - qs'/t + p\theta/t}} \, dv_{g} \nonumber\\
\leq \left(\int_{\Omega}\frac{\abs{u}^{qs'}}{\rho^{qs'}} \abs{\nabla \rho}^{qs'}\rho^{p\theta} dv_{g}\right)^{1/t} 
\left(\int_{\Omega} \frac{\abs{u}^{qs'}\abs{\nabla \rho}^{\left(\gamma + \epsilon\right)rs't' - qs't'/t}}{\rho^{\left(\gamma r+p^{*}\theta/s\right)s't' - qs't'/t + p\theta t'/t}} dv_{g}\right)^{1/t'}. \label{gen2}
\end{eqnarray}
Now, requiring that the following conditions are satisfied
\be qs' = p, \quad \left(r-q\right)s = p^{*}, \label{rrr}\ee
we get
\be s = \frac{p^{*}-p}{r-p} > 1, \label{rrr1}\ee
since $(\ref{condr})$ holds. 
Using $(\ref{rrr})$, by $(\ref{gen1})$ and $(\ref{gen2})$ we have 
\begin{eqnarray}
\int_{\Omega}\frac{\abs{u}^{r}}{\rho^{\gamma r}}\abs{\nabla \rho}^{\left(\gamma +\epsilon\right)r} dv_{g}  \leq 
\left(\int_{\Omega}\abs{u}^{p^{*}}\rho^{p^{*}\theta} dv_{g}\right)^{1/s} \nonumber\\
\left(\int_{\Omega}\frac{\abs{u}^{p}}{\rho^{p}} \abs{\nabla \rho}^{p}\rho^{p\theta} dv_{g}\right)^{1/s't} 
\left(\int_{\Omega} \frac{\abs{u}^{p}\abs{\nabla \rho}^{\left(\gamma + \epsilon\right)rs't' - pt'/t}}{\rho^{\left(\gamma r+p^{*}\theta/s\right)s't' - pt'/t + p\theta t'/t}} \, dv_{g}\right)^{1/s't'} \nonumber\\
\leq    C_2^{-p^*/s}   H\left(p\theta,p\right)^{-1/s't}\left(\int_{\Omega}\rho^{p\theta}\abs{\nabla u}^{p} dv_{g}\right)^{p^{*}/ps + 1/s't}   \nonumber\\
\left(\int_{\Omega} \frac{\abs{u}^{p}\abs{\nabla \rho}^{\left(\gamma + \epsilon\right)rs't' - pt'/t}}{\rho^{\left(\gamma r+p^{*}\theta/s\right)s't' - pt'/t + p\theta t'/t}} \, dv_{g}\right)^{1/s't'} ,
\end{eqnarray}
where, in the last inequality, we have used $(\ref{hsob})$ and and $(H_{\alpha})$ with $\alpha = p\theta$. 
To conclude we have to choose $t>1$ such  that 
\be   \left(\gamma + \epsilon\right)rs't' - pt'/t = p(\delta+\sigma) , \qquad 
\left(\gamma r+p^{*}\theta/s\right)s't' - pt'/t + p\theta t'/t = p\delta , \label{p1}\ee
and 
\be   p^{*}/ps + 1/s't = ar/p, \qquad    1/s't' = \left(1-a\right)r/p.     \label{p2}\ee 
First of all, note that in $(\ref{gen2})$ we can make the choice 
$t =\frac{sp}{s'\left(asr - p^{*}\right)} >1$, since $(\ref{condr})$ holds. 
Using the expressions of $t$ and $s$, equalities $(\ref{p2})$ 
follow by simple computations. 
Moreover, from $(\ref{p2})$ we obtain that 
$s't'=\frac{p}{r\left(1-a\right)}$ and $\frac{t'}{t}=\frac{asr -p^{*}}{sr\left(1-a\right)}$.
Using this two expressions and the conditions $(\ref{cond1})$ and $(\ref{cond2})$ 
we get also $(\ref{p1})$. 
This concludes the proof.
\hfill$\Box$

\boss
The condition $(\ref{cond2})$ takes into account the presence of the $\abs{\nabla\rho}$ 
in the weights appearing in $(\ref{disgen})$ and it is also a necessary condition. 
Indeed,  to see the necessity of $(\ref{cond2})$ we  argue as follows.
Assume that Theorem \ref{cknmanif} were true.
If $(S)$ and $(H_{\alpha})$ hold with a function $\rho$, then those inequalities still hold with the function $\lambda\rho$ for every $\lambda >0$, 
and hence the conclusion of Theorem \ref{cknmanif} holds
replacing $\rho$ with $\lambda\rho$. 
By homogeneous consideration one derives the necessity of  $(\ref{cond2})$.
\eoss

\boss
Since the condition $(\ref{cond2})$ is a requirement on the parameters $\epsilon$ and $\sigma$, 
if $\abs{\nabla \rho}=1$, these parameters do not appear in the inequality $(\ref{disgen})$. 
Therefore, condition $(\ref{cond2})$ is always fulfilled 
(i.e. choosing $\epsilon=a\theta$ and $\sigma=0$). 
The next corollary deals with a generalization of this case. 
\eoss

\bc
Assume that $(S)$ holds on $\Omega$ with $p^*>p$.
Let $\theta\in\Rset$ and $\rho\geq0$ be a function such that
$(H_{\alpha})$ holds with $\alpha=p\theta$ and  $\rho=d^\beta$ with
$\beta\in \Rset$ and $\abs{\nabla d}=1$.
Let $r>0$, $0\leq a\leq1$ and $\gamma$, $\delta$ be real numbers 
satisfying (\ref{condr}) and 
\be  \gamma + \frac{p^{*}\left(r-p\right)}{r\left(p^{*}-p\right)}=(1-\beta \theta) a+\delta\left(1-a\right).
\label{condiz1}\ee
Then there exists $C_3'>0$ such that  
\begin{eqnarray} C_3'\left(\int_{\Omega}\frac{\abs{u}^{r}}{d^{\gamma r}} dv_{g}\right)^{1/r} 
\leq  \left(\int_{\Omega}d^{\beta\theta p}\abs{\nabla u}^{p} \, dv_{g}\right)^{a/p}  
\left(\int_{\Omega}\frac{\abs{u}^{p}}{d ^{\delta p} }\, dv_{g}\right)^{\left(1-a\right)/p} , 
\quad u\in\Cinfinito_{0}\left(\Omega\right), 
\label{disge}\end{eqnarray}
that is 
\be C_3'\left|\frac{u}{\rho^{\gamma}}\right|_{L^{r}} \leq  \left|\rho^{\beta \theta}\abs{\nabla u}\right|^{a}_{L^{p}}\left|\frac{u}{\rho^{\delta}}\right|^{1-a}_{L^{p}}, \ee
\noindent provided $\frac{1}{\rho^{p\delta}}\in L^{1}_{loc}\left(\Omega\right)$. 
In particular, $C_3'=C_3 \beta^{\gamma+\beta\theta a-(1-a)\delta}$.
\ec

\boss
Notice that, if in the previous corollary we take $p^{*} = \frac{pN}{N-p}$, 
condition $(\ref{condiz1})$ becomes
\be \frac{1}{r} - \frac{\gamma}{N} = \frac{1}{p} + \frac{a}{N}\left(\beta\theta-1\right) - \frac{\delta}{N}\left(1-a\right), \label{cond3}\ee
and $(\ref{disge})$ is a 
particular case on manifold of a result obtained by Caffarelli, Kohn and Nirenberg in \cite{CKN} 
in Euclidean setting. 
\eoss

\section{Some Applications}\label{sec:ex}

In what follows we apply the results proved in the above sections to concrete cases, 
obtaining Hardy inequalities which in some cases are new. 
For sake of brevity we shall limit ourselves to show some applications of Theorems \ref{teo:g} and \ref{hardygen} 
by specializing the function $\rho$. 
With the same technique it is possible to obtain applications of the other theorems presented 
in the previous sections (Caccioppoli inequality, uncertain principle, Gagliardo-Nirenberg inequality, first order interpolation inequalities, and so on). 
We leave the details to the interested reader.

\subsection{Hardy inequality involving the distance from the boundary}

In order to prove a Hardy inequality involving the distance from the boundary, we need the following result, which is an immediate consequence of Theorems $\ref{teo:g}$ and $\ref{hardygen}$.
\bt 
Let $\overline{M}$ be a compact Riemannian manifold with boundary of class $\Cinfinito$, 
let $\varphi_{1}$ be the first eigenfunction related to the first eigenvalue of 
the $p$-Laplacian\footnote{This means that $(\lambda_1, \varphi_{1})$ is a solution of the problem 
\be  \left\{\begin{array}{ll} -\Delta_{p}\varphi = \lambda \abs{\varphi}^{p-2}\varphi &\mathrm{on\ } M, \cr \varphi= 0 &\mathrm{on\ } \partial M,\end{array}\right. 
\label{ep}\ee 
and $\lambda_1:=\min \{\lambda : (\lambda,\varphi)$ solves (\ref{ep}) $\}$.}, 
and $\alpha<p-1$. 
If $\varphi_{1}>0$, then the following inequality holds on $M$
\be \left(\frac{p-1-\alpha}{p}\right)^{p} \int_{M}\varphi_{1}^{\alpha} \frac{\abs{u}^{p}}{\varphi_{1}^{p}} \abs{\nabla\varphi_{1}}^{p}dv_{g} \leq  \int_{M}\varphi_{1}^{\alpha} \abs{\nabla u}^{p} dv_{g}, \quad u\in\Cinfinito_{0}\left(M\right). \ee

In particular we have
\be\label{dis:heigen} \left(\frac{p-1}{p}\right)^{p} \int_{M}\frac{\abs{u}^{p}}{\varphi_{1}^{p}} \abs{\nabla\varphi_{1}}^{p} dv_{g} \leq  \int_{M}\abs{\nabla u}^{p} dv_{g}, \quad u\in\Cinfinito_{0}\left(M\right). \ee
\et

The main theorem of this section is the following:
\bt\label{teo:hdist} Let $\overline{M}$ be a compact Riemannian manifold with boundary of class $\Cinfinito$, 
let $\varphi_{1}$ be the first eigenfunction related to the first eigenvalue of 
the $p$-Laplacian.
Assume that $\varphi_{1}\in\Cuno\left(\overline{M}\right)$, 
$\varphi_{1}>0$ on $M$ and $\abs{\nabla \varphi_{1}}\neq 0$ on $\partial M$.
 
Denoted by $d(x):=dist\left(x,\partial M\right)$, there exists a constant $c>0$ such that
\be c \int_{M}\frac{\abs{u}^{p}}{d^{p}}\, dv_{g} \leq  
\int_{M}\abs{\nabla u}^{p} dv_{g}, \quad u\in\Cinfinito_{0}\left(M\right). \ee
\et

The proof of the above theorem relies on the following result, which is worth of mention:
\bt 
Let $\overline{M}$ be a compact Riemannian manifold with boundary of class $\Cinfinito$, 
let $\varphi_{1}>0$ be the first eigenfunction related to the first eigenvalue  $\lambda_{1}$ of the $p$-Laplacian, and $0<s<p-1$. Then the following inequality holds
\be\label{dis:poinc} \lambda_{1}\frac{\left(p-1-s\right)^{\left(p-1\right)}}{p^{p}} \int_{M}\varphi_{1}^{s} \abs{u}^{p} dv_{g} \leq 
\int_{M} \varphi_{1}^{s} \abs{\nabla u}^{p} dv_{g}, \quad u\in\Cinfinito_{0}\left(M\right). \ee
\et
\noindent{\bf Proof.} 
Set $\phi:=\varphi_{1}^{s/\left(p-1\right)}$. By computation we have
$$-\Delta_{p} \phi= 
\frac{\abs{\nabla\phi}^{p}}{\phi}\left(p-1-s\right)\left(\frac{p-1}{s}\right) 
+ \lambda_1 \left(\frac{s}{p-1}\right)^{p-1} \phi^{p-1}.$$
Choosing $h:= -\abs{\nabla \phi}^{p-2}\nabla \phi$ and $A_h:= -\Delta_{p} \phi$, 
an application of Lemma $\ref{lemma:k1}$ yields
\bern \lambda_{1}\left(\frac{s}{p-1}\right)^{p-1} \int_M\phi^{p-1} \abs{u}^{p} dv_{g}\le \int_M (-\Delta_{p} \phi) \abs{u}^{p} dv_{g} \le  
p^{p} \int_M\frac{\abs{\nabla \phi}^{p\left(p-1\right)}}{\left( -\Delta_{p} \phi\right)^{p-1}}\abs{\nabla u}^{p} dv_{g}
\nonumber\\ \le p^{p} \left(\frac{s}{p-1}\right)^{p-1} \left(\frac{1}{p-1-s}\right)^{p-1}
\int_M\phi^{p-1} \abs{\nabla u}^{p} dv_{g}. 
\eern
This last chain of inequalities concludes the proof. 
\hfill $\Box$

\medskip

\noindent{\bf Proof of Theorem \ref{teo:hdist}.} For a fixed number $\gamma>0$,
we denote by $\Omega^{\gamma}$  and $\Omega_{\gamma}$ respectively the  sets
$\Omega^{\gamma}:=\{x : d\left(x\right) < \gamma\} $ and 
$\Omega_{\gamma}:=\{x : d\left(x\right)\ge \gamma\} $.

Let $\varphi_1$ be such that $\abs{\abs{\varphi_1}}_\infty\le 1.$
By continuity argument, we have that there exist $\epsilon,b>0$ such that
$$ \abs{\nabla\varphi_{1}\left(x\right)}\ge b>0,\qquad\mathrm{for}\  x\in\Omega^{\epsilon}.$$

Since $\varphi_{1}$ is a Lipschitz continuous function, we obtain that there exist $L>0$ such that
$$ \varphi_{1}\left(x\right)\le L d\left(x\right),\qquad x\in M.$$

We set
$$l_{\epsilon}:=\min_{\Omega_{\epsilon}} \varphi_{1}.$$
From $(\ref{dis:heigen})$ and $(\ref{dis:poinc})$ we get respectively
$$ \int_{M}\abs{\nabla u}^{p} dv_{g} \ge 
\left(\frac{p-1}{p}\right)^{p} \int_{\Omega^\epsilon}\frac{\abs{u}^{p}}{\varphi_{1}^{p}}\abs{\nabla\varphi_{1}}^{p}dv_{g}
\ge \left(\frac{p-1}{p}\right)^{p} \frac{b^{p}}{L^{p}}\int_{\Omega^{\epsilon}}\frac{\abs{u}^{p}}{d^{p}}\, dv_{g} ,  $$
$$ \int_{M}\abs{\nabla u}^{p} dv_{g} \ge 
\lambda_{1}\frac{\left(p-1-s\right)^{\left(p-1\right)}}{p^{p}}\int_{\Omega_{\epsilon}}\varphi_{1}^{s} \abs{u}^{p} dv_{g} \ge \lambda_{1}\frac{\left(p-1-s\right)^{\left(p-1\right)}}{p^{p}}
l_{\epsilon}^{s} \epsilon^{p} \int_{\Omega_{\epsilon}} \frac{\abs{u}^{p}}{d^{p}}\, dv_{g}.
$$
Choosing $2c:=\min\left\lbrace  \left(\frac{p-1}{p}\right)^{p} \frac{b^{p}}{L^{p}} ,\lambda_{1}\frac{\left(p-1-s\right)^{\left(p-1\right)}}{p^{p}}
\, l_{\epsilon}^{s} \epsilon^{p} \right\rbrace$ and summing up the above estimates we obtain the claim.
\hfill$\Box$

\subsection{Hardy inequality for $p$-hyperbolic manifold}

In this section we establish Hardy inequalities involving the Green function of the operator $-\Delta_{p}$. 
The case $p=2$ is already proved in \cite{LW}. 

Examples of $p$-hyperbolic manifolds are the following. 
The Euclidean space $\Rset^N$ is $p$-hyperbolic for $N>p$.
If $M$ is a Cartan-Hadamard manifold (see Section \ref{sec:Cartan})
whose sectional curvature $K_M$ is uniformly negative, that is $K_M\le -a^2<0$,
then $M$ is $p$-hyperbolic for any $p>1$ (see \cite{Ho} and \cite{HP}).

We have the following.
\bt\label{teo-par} 
Let $\left(M, g\right)$ be a $p$-hyperbolic manifold, 
let $G_x$ be the Green function for $\Delta_p$ with pole at $x$, and $\alpha\in\Rset$. 
Then the following inequality holds 
\be \left(\frac{\abs{p-1-\alpha}}{p}\right)^{p} \int_{M\setminus\left\{x\right\}} G_{x}^{\alpha}\frac{\abs{\nabla G_{x}}^{p}}{G_{x}^{p}} \abs{u}^{p} dv_{g} \leq \int_{M\setminus\left\{x\right\}}G_{x}^{\alpha}\abs{\nabla u}^{p} dv_{g}, \quad u\in\Cinfinito_{0}\left(M\setminus\left\{x\right\}\right). \label{Hpar1}\ee
In particular, we have
\be \left(\frac{p-1}{p}\right)^{p} \int_{M\setminus\left\{x\right\}}\frac{\abs{\nabla G_{x}}^{p}}{G_{x}^{p}} \abs{u}^{p} dv_{g} \leq \int_{M\setminus\left\{x\right\}}\abs{\nabla u}^{p} dv_{g}, \quad u\in\Cinfinito_{0}\left(M\setminus\left\{x\right\}\right). \label{Hpar}\ee
and, if $p<N$, 
\be \left(\frac{p-1}{p}\right)^{p}
\int_{M} \frac{\abs{\nabla G_{x}}^{p}}{G_{x}^{p}} \abs{u}^{p} dv_{g} \leq\int_{M}\abs{\nabla u}^{p} dv_{g}, \quad u\in\Cinfinito_{0}\left(M\right). \label{Hpar2}\ee
\et

\noindent{\bf Proof.} 
We know that $G_{x}\in W^{1, p}_{loc}\left(M\right)$ 
is a nonnegative function on $M$.
Moreover the hypotheses of Theorem $\ref{hardygen}$ are fulfilled; in fact
\begin{enumerate} 
\item[i)]$-\Delta_{p}G_{x}=0$ in $M\setminus\left\{x\right\}$
\item[ii)]$\frac{\abs{\nabla G_{x} }^{p}}{G_{x}^{p-\alpha}}, G_{x}^{\alpha}\in L^{1}_{loc}\left(M\setminus\left\{x\right\}\right).$
\end{enumerate}
Then, by Theorem $\ref{hardygen}$, inequality $(\ref{Hpar1})$ holds.
In particular, taking $\alpha =0$, we obtain the inequality $(\ref{Hpar})$. 
Moreover, if $p<N$, we are in the position to apply Corollary $\ref{coro:q1}$, 
because $\left\{x\right\}$ is a set of zero $p$-capacity 
(see Theorem 2.27 in \cite{HKM}), 
and then we can use Proposition \ref{appe}; this proves that also $(\ref{Hpar2})$ holds.  
\hfill$\Box$

\subsection{Hardy inequality on Cartan-Hadamard manifold}\label{sec:Cartan}

In what follows $\left(M, g\right)$ will denote a Cartan-Hadamard manifold, that is, a connected, simply connected, complete 
Riemannian manifold of dimension $N\geq2$, of nonpositive sectional curvature 
(see \cite{Ca, Ho, LW} for further details). 
Examples of Cartan-Hadamard manifolds are the Euclidean space~$\Rset^N$ with the usual metric 
(which has constant sectional curvature equal to zero), and the standard $N$-dimensional 
hyperbolic space~$\Hset^N$ (which has constant sectional curvature equal to~$-1$). 
 
Let $o\in M$ be a fixed point and denote by $r$ the distance function from $o$. 
We have the following.
\bt\label{teo:h}
Let $\left(M, g\right)$ be a Cartan-Hadamard manifold and $\alpha\in\Rset$. 

If $(N-p)(p-1-\alpha)>0$, we have
\begin{eqnarray}
\left(\frac{\left(N-p\right)\left(p-1-\alpha\right)}{p\left(p-1\right)}\right)^{p} \int_{M\setminus\left\{o\right\}} r^{\alpha\frac{p-N}{p-1}} \frac{\abs{u}^{p}}{r^{p}} \, dv_{g} \leq \int_{M\setminus\left\{o\right\}}r^{\alpha\frac{p-N}{p-1}} \abs{\nabla u}^{p} dv_{g}, \nonumber\\
u\in\Cinfinito_{0}\left(M\setminus\left\{o\right\}\right). \label{cartan1}\end{eqnarray}

If $1<p<N$, we have
\be\left(\frac{N-p}{p}\right)^{p} \int_{M}\frac{\abs{u}^{p}}{r^{p}}\, dv_{g} \leq \int_{M}\abs{\nabla u}^{p} dv_{g}, \quad
u\in\Cinfinito_{0}\left(M\right). \label{cartan3}\ee

If $1<p\le N$, setting $\Omega\decl r^{-1}(\left[0,1\right[)$, we have
\be\left(\frac{{p-1}}{p}\right)^{p} \int_{\Omega} \frac{\abs{u}^{p}}{\abs{r\ln r}^{p}}\, dv_{g} \leq \int_{\Omega} \abs{\nabla u}^{p} dv_{g}, \quad
u\in\Cinfinito_{0}\left(\Omega\right). \label{cartan3log}\ee

Let  $1<p\le N$. If $p-1>\alpha$ we set $\Omega\decl r^{-1}(\left]0,1\right[)$, 
else if $p-1<\alpha$ we set  $\Omega\decl r^{-1}(\left]1, +\infty\right[)$. 
We have
\be\left(\frac{\abs{p-1-\alpha}}{p}\right)^{p} \int_{\Omega}\abs{\ln r}^\alpha \frac{\abs{u}^{p}}{\abs{r\ln r}^{p}}\, dv_{g} 
\leq \int_{\Omega}\abs{\ln r}^\alpha\abs{\nabla u}^{p} dv_{g}, \quad
u\in\Cinfinito_{0}\left(\Omega\right). \label{cartan4log}\ee
\et

\noindent{\bf Proof.} 
In $M\setminus\left\{o\right\}$ we define $\rho =r^{\beta}$ with $\beta\in\Rset$ that will be chosen later. 
The function $\rho \in W^{1, p}_{loc}\left(M\setminus\left\{o\right\}\right)$ is nonnegative on $M\setminus\left\{o\right\}$. 

The function $r$ satisfies the relations
\be \abs{\nabla r}=1, \qquad \Delta r\geq\frac{N-1}{r} \label{cartan}\ee
(see \cite{Ca}). 

By computation we obtain
\begin{eqnarray}
\Delta_{p} \rho &=& \Div\left(\abs{\nabla\rho}^{p-2}\nabla\rho\right) \nonumber\\
&=& \Div\left(\abs{\beta}^{p-2} r^{\left(\beta -1\right)\left(p-2\right)} \beta r^{\beta -1} \nabla r\right)\nonumber\\
&=& \abs{\beta}^{p-2} \beta  \Div\left(r^{\left(\beta -1\right)\left(p-1\right)} \nabla r \right) \nonumber\\
&=& \abs{\beta}^{p-2} \beta \left[\left(\beta -1\right)\left(p-1\right) r^{\left(\beta -1\right)\left( p-1\right) -1}
+ r^{\left(\beta -1\right)\left( p-1\right)} \Delta r\right] \nonumber\\
&=& \abs{\beta}^{p-2} \beta r^{\left(\beta -1\right)\left( p-1\right) -1} 
\left[\left(\beta -1\right)\left( p-1\right) + r\Delta r\right], \nonumber\end{eqnarray}
and hence
\be 
-\left( p-1-\alpha\right) \Delta_{p} \rho = 
-\left( p-1-\alpha\right) \abs{\beta}^{p-2}\beta r^{\left(\beta -1\right)\left( p-1\right) -1} 
\left[ \left(\beta -1\right)\left( p-1\right) + r\Delta r\right]. \label{cart}\ee

Next choosing $\beta =\frac{p-N}{p-1}$, using (\ref{cartan}), we have that 
$\left(\beta -1\right)\left( p-1\right) + r\Delta r\geq 0$, 
and then $-\left( p-1-\alpha\right) \Delta_{p} \rho\geq 0$. That is 
the hypothesis i) of Theorem $\ref{hardygen}$ is fulfilled. 
Since $\frac{\abs{\nabla \rho}^{p}}{\rho^{p-\alpha}} = \abs{\beta}^{p} \frac{1}{r^{p-\alpha}}\in L^{1}_{loc}\left(M\setminus\left\{o\right\}\right)$ and 
$\rho^{\alpha}=r^{\beta\alpha}\in L^{1}_{loc}\left(M\setminus\left\{o\right\}\right)$, 
we are in a position to apply Theorem \ref{hardygen}, 
obtaining the inequality $(\ref{cartan1})$. 
 
In particular, taking $\alpha =0$ and using Corollary $\ref{coro:q1}$ we get $(\ref{cartan3})$. 
Indeed, since $p\le N$, $\left\{o\right\}$ is a set of zero $p$-capacity 
(see Theorem 2.27 in \cite{HKM}), 
and then we can use Proposition \ref{appe}. 

Next we prove (\ref{cartan4log}). To this end, by choosing
$\rho\decl (\alpha-p+1)\ln r$, we have that $\rho>0$ in $\Omega$ (according to the different
cases $\alpha>(<) p-1$). By computation we have
\begin{eqnarray*}
-(p-1-\alpha) \Delta_{p} \rho&=& (\alpha-p+1)\ \Div\left(\abs{\nabla\rho}^{p-2}\nabla\rho\right) \\ 
&=& (\alpha-p+1)^p \ \Div( r^{1-p}\nabla r)\\
&=& (\alpha-p+1)^p \ \left( \frac{1-p}{r^p}+ \frac{r\Delta r}{r^p} \right)\\
&\ge& (\alpha-p+1)^p \left( \frac{N-p}{r^p}\right) \ge 0.
\end{eqnarray*}
The claim follows applying Theorem $\ref{hardygen}$.

We conclude the proof by proving (\ref{cartan3log}).
Choosing $\alpha=0$ in (\ref{cartan4log}), 
we have that inequality $(\ref{cartan3log})$ holds 
for every $\in\Cinfinito_{0}\left(\Omega\setminus \left\lbrace o\right\rbrace \right)$.
However in this case 
$\left\lbrace o\right\rbrace$ 
is a set of zero $p$-capacity and, applying Corollary $\ref{coro:q1}$, we complete the proof.
\hfill$\Box$
\medskip

Inequality (\ref{cartan3}) is present in \cite{Ca} for $p=2$ .
In \cite{KO} the authors prove (\ref{cartan1}) for $p=2$ and for a special case
of manifold $M$, namely, when $M$ is the unit ball modeling the  standard hyperbolic space $\Hset^N$.
For this case the authors prove that the constant in (\ref{cartan1}) is sharp
and they show that  a remainder term can be added.

\subsection{Hardy inequalities involving the distance from the soul of a manifold}

Let $\left(M, g\right)$ be a complete non compact Riemannian manifold, 
of dimension $N\geq2$, with nonnegative sectional curvatures. 
A result due to Cheeger and Gromoll asserts that 
there exists a compact embedded totally convex submanifold $S$ 
with empty boundary, whose normal bundle is diffeomorphic to $M$ (see \cite{CG}). 
The submanifold $S$, called "`soul"' of $M$, 
is not necessarily unique but every two souls of $M$ are isometric. 
\emph{``Totally convex''} means that any geodesic arc in $M$ connecting 
two points in $S$ (which may coincide) lies entirely in $S$. 
In particular, $S$ is connected, totally geodesic in $M$, 
and has nonnegative sectional curvature. 
Moreover $0\leq dim S<dim M$.

Denote by $r : M\setminus S\rightarrow\Rset$ the distance function to $S$. 
We have that $r$ is smooth on $M\setminus S$ and $\abs{\nabla r} = 1$ on $M\setminus S$. 
Now we suppose that radial sectional curvature $K_r$, that is 
sectional curvature of two-planes containing the direction $\nabla r$, 
satisfies  
\be 0\leq K_{r}\leq \frac{c_{N}\left(1-c_{N}\right)}{r^{2}}, \label{curv}\ee
where $c_{N} = \frac{N-2}{N}$; 
then we have
\be \Delta r\geq \frac{c_{N}\left(N-s-1\right)}{r}, \label{curv1}\ee
where $s=dimS$ (see \cite{EF}). We have the following:
\bt 
Let $\left(M, g\right)$ be a Riemannian manifold with 
nonnegative curvature. 
Suppose that $(\ref{curv})$ is fulfilled. 

Let $G\decl c_{N}\left(N-s-1\right)-p+1$. 
If $G\cdot (p-1-\alpha)>0$, we have
\begin{eqnarray}
\left(\frac{G\cdot \left(p-1-\alpha\right)}{p\left(p-1\right)}\right)^{p} 
\int_{M\setminus S}r^{-\alpha\frac{G}{p-1}} \frac{\abs{u}^{p}}{r^{p}} \, dv_{g} \leq
\int_{M\setminus S}r^{-\alpha\frac{G}{p-1}}\abs{\nabla u}^{p} dv_{g}, \nonumber\\
\quad u\in\Cinfinito_{0}\left(M\setminus S\right). \label{cn}\end{eqnarray}
Moreover, if $G>0$, we have
\be \left(\frac{G}{p}\right)^{p} \int_{M} \frac{\abs{u}^{p}}{r^{p}}\, dv_{g} \leq \int_{M}\abs{\nabla u}^{p} dv_{g}, \quad u\in\Cinfinito_{0}\left(M\right). \label{cn1}\ee
\et

\noindent{\bf Proof.} 
Let $\rho=r^{\beta}$ in $M\setminus S$ with $\beta =-\frac{G}{p-1}$. 
Arguing as in the proof of Theorem \ref{teo:h}, 
using $(\ref{cart})$ and $(\ref{curv1})$, 
we obtain $-\left( p-1-\alpha\right) \Delta_{p} \rho\geq 0$. 
Since $\rho\in W^{1,p}_{loc}(M\setminus S)$, 
$\frac{\abs{\nabla \rho}^{p}}{\rho^{p-\alpha}} = \abs{\beta}^{p} \frac{1}{r^{p-\alpha}}\in L^{1}_{loc}\left(M\setminus S\right)$ and $\rho^{\alpha}=r^{\beta\alpha}\in L^{1}_{loc}\left(M\setminus S\right)$, 
the hypotheses of Theorem \ref{hardygen} are fulfilled, 
and $(\ref{cn})$ follows.

In particular, if $\alpha =0$, $(\ref{cn})$ becomes 
\be \left(\frac{G}{p}\right)^{p} \int_{M\setminus S}\frac{\abs{u}^{p}}{r^{p}} \, dv_{g}\leq \int_{M\setminus S}\abs{\nabla u}^{p} dv_{g}, \quad u\in\Cinfinito_{0}\left(M\setminus S\right). \nonumber\ee
Now, the hypothesis $G>0$ implies that $N-p>s$. 
In fact, by the fact that $G = c_N\left(N-s-1\right) -p+1 = \frac{N-2}{N} \left(N-s-1\right) -p+1 >0$, 
by simple computations we get
\bern
\left(N-2\right)\left(N-s-1\right) -N\left( p+1\right) &=& 
\left(N-2\right)\left(N-s\right) -Np   \nonumber\\
&=& N \left(N-s\right) -2\left(N-s\right) -N p \nonumber\\
&=& N\left(N-s-p\right)-2\left(N-s\right) >0,
\nonumber  \eern
which implies $N-s-p > 2\frac{N-s}{N}>0$. 
Then $S$ is a set of zero $p$-capacity 
(see Theorem 2.27 in \cite{HKM}), 
and we can use Proposition \ref{appe} and Corollary $\ref{coro:q1}$ to obtain inequality $(\ref{cn1})$. 
\hfill$\Box$

\subsection{Hardy-Poincar\'e inequality for the hyperbolic plane} 
Let $\Cset_{+}=\left\{z=x+iy : Im z=y>0\right\}$ be  the upper half-plane
equipped with the Poincar\'{e} metric $ds^{2}=\frac{dx^{2}+dy^{2}}{y^{2}}$. 
This space is a Riemannian manifold modeling the \emph{two dimensional hyperbolic space}. 
In this case, the gradient~$\nabla_H$, the divergence~$\Div_{\!\!H}$, 
the Laplacian~$\Delta_H$ and the volume $dv_g$ related to the metric 
are respectively the following
\begin{equation}\begin{split}\label{intrinsic}
\nabla_H u &= y\nabla_E u, \\
\Div_{\!\!H} &= y^2 \Div_{\!\!E} , \\
\Delta_H u &= y^2 \Delta_E u, \\
dv_g&=\frac{dx\, dy}{y^2},
\end{split}\end{equation}
where we have denoted with $\nabla_E$, $\Div_{\!\!E}$, $\Delta_E$ the related operator in the Euclidean setting, and $dx\,dy$ is the Lebesgue measure in $\Rset^2$. 

By using Theorem $\ref{hardygen}$ with $p=2$, we deduce a 
Hardy inequality on the upper half-plane.
\bt
Let $\alpha\in\Rset$. For every $u\in\Cinfinito_{0}\left(\Cset_{+}\right)$ we have
$$ \frac{(1-\alpha)^2}{4}\int_{\Cset_{+}} y^{\alpha}|u|^2 \frac{dx\, dy}{y^2}\leq \int_{\Cset_{+}} y^{\alpha}|\nabla_H u|^2 \frac{dx\, dy}{y^2}, \quad u\in\Cinfinito_{0}\left(\Cset_{+}\right). $$
\et

\noindent{\bf Proof.} We consider the function $\rho(z)=y$, where $z=x+iy$. 
Clearly, $\rho$ belongs to $W^{1,2}_{loc}(\Cset_{+})$, and~$\rho^{\alpha}=y^{\alpha}$ 
belongs to~$L^1_{loc}(\Cset_{+})$. 
Moreover, from~\eqref{intrinsic}, we have that 
$$ \frac{|\nabla_H \rho|^2}{\rho^{2-\alpha}}=\frac{y^2|\nabla_E\rho|^2}{y^{2-\alpha}}=y^{\alpha}\in L^1_{loc}(\Cset_{+}), $$
and 
$$ \Delta_H \rho = y^2 \Delta_E\rho =0. $$ 
Therefore, the hypotheses of Theorem~\ref{hardygen} are satisfied and this 
concludes the proof. 
\hfill$\Box$

\subsection{The Euclidean case}\label{sec:euclidean}
In this last section we show that our main results, Theorems $\ref{teo:g}$ and $\ref{hardygen}$, yield some well known sharp Hardy inequalities in the Euclidean space.  

\bigskip 

Since $\Rset^N$ is $p$-hyperbolic for $N>p$ and it is a Cartan-Hadamard manifold, 
Theorems \ref{teo-par} and \ref{teo:h} hold also on $\Rset^N$ with 
$G_x = \vert x\vert^{\frac{p-N}{p-1}}$ and $r=\vert x\vert$ respectively. 
However, the function $\vert x\vert^{\frac{p-N}{p-1}}$ is $p$-harmonic in $\Rset^N\setminus\left\lbrace 0\right\rbrace$ also for $p>N$. 
Therefore we have the following:
\bt\label{teo:h2} 
Let $p\neq N$. For any $u\in\Cinfinito_{0}\left(\Rset^{N}\setminus\left\{0\right\}\right)$ we have
\bern
\left(\frac{\abs{N-p}\cdot \abs{p-1-\alpha}}{p\left(p-1\right)}\right)^{p} 
\int_{\Rset^{N}\setminus\left\{0\right\}}    \abs{x}^{\alpha\frac{p-n}{p-1}}
\frac{\abs{u}^{p}}{\abs{x}^{p}} \,dx 
\leq \int_{\Rset^{N}\setminus\left\{0\right\}}   \abs{x}^{\alpha\frac{p-n}{p-1}} \abs{\nabla u}^{p} dx.
\label{omega}\eern
\et

\bigskip

In the half space $\Rset^{N}_{+}$ there holds the following:
\bt\label{Thalf} 
Let $\alpha\in\Rset$, let $N\geq 2$, 
let $\Rset^{N}_{+} = \left\{\left(x_{1}, \ldots, x_{N}\right)\in\Rset^{N} : x_{1}>0\right\}$, 
and let $\rho\left(x\right)\decl d\left(x, \partial \Rset^{N}_{+}\right)$ 
be the distance from the boundary of $\Rset^{N}_{+}$. 
Then we have 
\be \left(\frac{\abs{p-1-\alpha}}{p}\right)^{p} \int_{\Rset^{N}_{+}} \rho^{\alpha}\frac{\abs{u}^{p}}{\rho^{p}} \, dx\leq  \int_{\Rset^{N}_{+}}\rho^{\alpha}\abs{\nabla u}^{p} dx, 
\quad u\in\Cinfinito_{0}\left(\Rset^{N}_{+}\right). \label{half}\ee
\et

\noindent{\bf Proof.} 
The distance $\rho\left(x\right)= x_1 \in W^{1, p}_{loc}\left(\Rset^{N}_{+}\right)$ is nonnegative on 
$\Rset^{N}_{+}$ and it is easy to verify that the hypotheses of Theorem $\ref{hardygen}$ 
are satisfied. Therefore the thesis follows.
\hfill$\Box$

\medskip

From Theorem~\ref{Thalf}, we can deduce, as a particular case, 
a well-known Hardy inequality for the upper 
half-plane $\Cset_{+}=\left\{z=x+iy : Im z=y>0\right\}$ 
(see for instance \cite{He} and references therein). 
\bc 
Let $\alpha\in\Rset$. Then, for every $u\in\Cinfinito_{0}\left(\Cset_{+}\right)$, we have the following:
\be
\left(\frac{\vert p-1-\alpha\vert}{p}\right)^{p} \int_{\Cset_{+}}\abs{u}^{p}   \frac{dA\left(z\right)}{\left(Im z\right)^{p-\alpha}}  \leq  2^{p/2}  \int_{\Cset_{+}} \left(Im z\right)^{\alpha} \left( \abs{\partial u\left(z\right)}^{2} + 
\abs{\overline{\partial} u\left(z\right)}^{2}\right)^{p/2} dA\left(z\right),      
\label{HP}\ee
where $dA\left(z\right):=\frac{dxdy}{\pi}$, and 
$\partial$, $\overline{\partial}$ are the Wirtinger operators, that is
\be      \partial := \frac{1}{2} \left( \frac{\partial}{\partial x} - i \frac{\partial}{\partial y}\right),   \qquad 
\overline{\partial} := \frac{1}{2} \left( \frac{\partial}{\partial x} + i \frac{\partial}{\partial y}\right). 
\label{de-debar}\ee
\ec

\bigskip

Actually a more general theorem for convex domains holds.  
\bt  
Let $D$ be a proper convex open subset of $\Rset^{N}$, let $d:=dist\left(\cdot , \partial D\right)$ be the distance from $\partial D$ and let $\alpha<p-1$. 
Then we have
\be \left(\frac{p-1-\alpha}{p}\right)^{p} \int_{D}d^{\alpha} \frac{\abs{u}^{p}}{d^{p}} \, dx\leq \int_{D}d^{\alpha} \abs{\nabla u}^{p} dx, \quad u\in\Cinfinito_{0}\left(D\right), \nonumber\ee 
and, in particular,
\be \left(\frac{p-1}{p}\right)^{p} \int_{D}\frac{\abs{u}^{p}}{d^{p}} \, dx\leq \int_{D}\abs{\nabla u}^{p} dx, \quad u\in\Cinfinito_{0}\left(D\right). \nonumber\ee 
\et

\noindent{\bf Proof.} The thesis will follow by applying Theorems \ref{teo:g} and \ref{hardygen}. To this end it suffices to prove that the function $d\left(x\right) :=dist\left(x, \partial D\right)$ is $p$-superharmonic. Indeed, $D=\bigcap\Pi$, and $d\left(x\right)=\inf_{\Pi} dist\left(x,\partial\Pi\right)$ where the intersection and the infimum are taken over all the half-spaces $\Pi$ containing $D$.
Since $dist\left(x,\partial\Pi\right)$ is continuous
and $p$-harmonic, we have that $d$ is $p$-superharmonic (see \cite{HKM}). 
This concludes the proof.
\hfill$\Box$

\appendice{Appendix}{A}

Let us recall that $p$-capacity of a compact set $K$ is defined as 
\begin{eqnarray}
cap_{p}\left(K, M\right) = \inf\left\{\int_{M}\abs{\nabla u}^{p} dv_{g} : u\in\Cinfinito_{0}\left(M\right),\ 0\leq u\leq1,\right.\nonumber\\ 
\left. u=1\ \mathrm{in\ a\ neighborhood\ of\ K\ }\right\}. \label{cap}\end{eqnarray}

\bpr\label{appe}
Let $M$ be a $p$-hyperbolic manifold of dimension $N$.  
Let $K\subset M$ be a compact set of zero $p$-capacity. Then 
\be D^{1,p}\left(M\right) \subset D^{1,p}\left(M\setminus K\right), \nonumber\ee
that is every function $u\in D^{1,p}\left(M\right)$ can 
be approximated by function $\Cinfinito_{0}\left(M\setminus K\right)$ 
in the norm $\abs{\cdot}_{D^{1,p}}$.
\epr

\noindent{\bf Proof.} 
Let $\varphi\in\Cinfinito_{0}\left(M\right)$. 
In order to prove the claim it is sufficient to prove that 
$\varphi\in D^{1,p}\left(M\setminus K\right)$. 
Since $cap_{p}\left(K, M\right)=0$, 
there exists a sequence $\left(u_{j}\right)_{j\geq1}$ such that, for any $j\geq1$, 
$u_{j}\in\Cinfinito_{0}\left(M\right)$, $0\leq u_{j}\leq1$,
$u_{j}=1$ in a neighborhood of $K$ and $u_{j}\rightarrow 0$ in $D^{1,p}\left(M\right)$. 
For every $j\geq1$ the function $\varphi_{j} := \left(1-u_{j}\right)\varphi$ 
belongs to $\Cinfinito_{0}\left(M\setminus K\right) \subset D^{1,p}\left(M\setminus K\right)$. 
We shall prove that  $\varphi_{j}\rightarrow\varphi$ in $D^{1,p}\left(M\setminus K\right)$, that is 
\be \int_{M\setminus K} \abs{\nabla\varphi_{j} - \nabla\varphi}^{p} dv_{g} \longrightarrow 0, \qquad  \mathrm{(as\ } j\rightarrow+\infty). 
\label{u}\ee 
In fact, we have  
\begin{eqnarray} 
\left(\int_{M\setminus K} \abs{\nabla\varphi_{j} - \nabla\varphi}^{p} dv_{g}\right)^{1/p} 
= \left(\int_{M\setminus K} \abs{\nabla\varphi\left(1-u_{j}\right) - \varphi\nabla u_{j} - \nabla\varphi}^{p} dv_{g}\right)^{1/p} \nonumber\\
\leq \left(\int_{M\setminus K} \abs{\nabla\varphi}^{p} \abs{u_{j}}^{p} dv_{g}\right)^{1/p} + 
\left(\int_{M\setminus K} \abs{\varphi}^{p} \abs{\nabla u_{j}}^{p} dv_{g}\right)^{1/p}. \label{uj}\end{eqnarray}
The second term in (\ref{uj}) converges to $0$ for $j\rightarrow+\infty$. 
Indeed, since $u_{j}\rightarrow 0$ in $D^{1,p}\left(M\right)$, we obtain
\bern
\int_{M\setminus K} \abs{\varphi}^{p} \abs{\nabla u_{j}}^{p} dv_{g} &\leq& 
\int_{M} \abs{\varphi}^{p} \abs{\nabla u_{j}}^{p} dv_{g} \\
&\leq& \left|\varphi\right|^{p}_{\infty} \int_{M} \abs{\nabla u_{j}}^{p} dv_{g} \longrightarrow 0, \qquad \mathrm{(as\ } j\rightarrow+\infty). 
\eern 
It remains 
to prove that the first term in (\ref{uj}) converges to $0$ as well. 
Let $D$ be the support of $\varphi$; then we get
\begin{eqnarray} 
\int_{M\setminus K} \abs{\nabla\varphi}^{p} \abs{u_{j}}^{p} dv_{g} \leq 
\int_{M} \abs{\nabla\varphi}^{p} \abs{u_{j}}^{p} dv_{g} 
= \int_{D} \abs{\nabla\varphi}^{p} \abs{u_{j}}^{p} dv_{g} \nonumber\\
\leq \left|\nabla\varphi\right|^{p}_{\infty} \int_{D} \abs{u_{j}}^{p} dv_{g} 
\leq \left|\nabla\varphi\right|^{p}_{\infty} C \int_{M}\abs{\nabla u_{j}}^{p} dv_{g} \longrightarrow 0, \qquad \mathrm{(as\ } j\rightarrow+\infty), 
\nonumber\end{eqnarray}
where, in the last inequality, we have used 
a characterization of the $p$-hyperbolic manifold 
(see Theorem $3$ in \cite{Tr2}). 
\hfill$\Box$

\end{document}